\documentclass[dvips,preprint,authoryear,12pt]{imsart}

\RequirePackage{natbib,amssymb,verbatim,undertilde,mathrsfs}
\RequirePackage[OT1]{fontenc}
\RequirePackage{amsthm,amsmath,natbib}
\RequirePackage{hypernat}
\usepackage{graphicx}
\def\Th{\Theta}

\def\Cov{\mbox{Cov}}

\def\defn{\triangleq}
\def\a{\alpha}

\def\dd{\boldsymbol{\delta}}

\def\th{\theta}

\def\tr{\mathrm{tr}}

\def\1{\mathbf{1}}
\def\aa{\boldsymbol{\alpha}}

\def\l{\lambda}

\def\e{\epsilon}
\def\ee{\boldsymbol{\epsilon}}
\def\b{\beta}

\def\bb{\boldsymbol{\beta}}
\def\d{\delta}

\def\x{\mathbf{x}}
\def\w{\mathbf{w}}
\def\F{\mathbb{F}}

\def\s{\sigma}

\def\u{\mathbf{u}}
\def\y{\mathbf{y}}

\def\1{\mathbf{1}}

\def\R{\mathbb{R}}

\def\X{\mathbf{X}}
\def\S{\mathit{\Sigma}}

\begin{document}

\bibliographystyle{ims}
\begin{frontmatter}
\title{Dense Signals, Linear Estimators, and Out-of-Sample Prediction for High-Dimensional Linear Models \protect}
\runtitle{Linear Estimators and High-Dimensional Linear Models}

\begin{aug}
\author{\fnms{Lee} \snm{Dicker}
\ead[label=e1]{ldicker@stat.rutgers.edu}}


\runauthor{L. Dicker}

\affiliation{Rutgers University}

\address{Department of Statistics and Biostatistics \\ Rutgers University \\ 501 Hill Center, 
 110 Frelinghuysen Road \\ Piscataway, NJ 08854 \\
\printead{e1}}
\end{aug}

\begin{keyword}[class=AMS]
\kwd[Primary ]{62J07}
\kwd[; secondary ]{62C20}
\kwd{15B52}
\end{keyword}

\begin{keyword}
\kwd{ridge regression}
\kwd{James-Stein estimator}  
\kwd{marginal regression}
\kwd{random matrix theory}
\kwd{oracle estimators}
\kwd{minimax estimators}
\end{keyword}

\begin{abstract}
Motivated by questions about dense (non-sparse) signals in high-dimensional data
analysis, we study the unconditional out-of-sample prediction error (predictive risk) associated with three popular linear estimators for high-dimensional linear models: ridge regression estimators, scalar multiples of the ordinary least squares (OLS) estimator (referred to as James-Stein shrinkage estimators), and marginal regression estimators.  The results in this paper require no assumptions about sparsity and imply: (i) if prior information about the population predictor covariance is available, then the ridge estimator outperforms the OLS, James-Stein, and marginal estimators; (ii) if little is known about the population predictor covariance, then the James-Stein estimator may be an effective alternative to the ridge estimator; and (iii) the marginal estimator has serious deficiencies for out-of-sample prediction.  
Both finite sample and asymptotic properties of the estimators are studied in this paper.  Though various asymptotic regimes are considered, we focus on the setting where the number of predictors is roughly proportional to the number of observations.    Ultimately, the results presented here provide new and detailed practical guidance regarding several well-known non-sparse methods for high-dimensional linear models.
\end{abstract}

\end{frontmatter}

{\allowdisplaybreaks
\section{Introduction}

High-dimensional data analysis is one of the most active areas of current statistical research. Much of this research has been driven by technological advances across a variety of scientific disciplines, including molecular biology and genomics, that have enabled investigators to collect vast datasets with relative ease.  The linear model has played a prominent role in recent literature on high-dimensional data analysis.  In the linear model, observed outcomes $y_1,...,y_n \in \R$ and corresponding $d$-dimensional predictors $\x_1,...,\x_n \in \R^d$ are related via the equation
\begin{equation} \label{lm}
y_i = \x_i^T\bb + \e_i, \ 1 \leq i \leq n,
\end{equation}
where $\bb = (\b_1,...,\b_d)^T \in \R^d$ is an unknown parameter vector, and $\e_1,...,\e_n$ are unobserved iid error terms with mean 0 and variance $\s^2 > 0$.  To simplify notation, let $\y = (y_1,...,y_n)^T \in \R^n$, $X = (\x_1,...,\x_n)^T$, and $\ee = (\e_1,...,\e_n)^T$.  Then the observed data are $(\y,X)$ and (\ref{lm}) may be rewritten as $\y = X\bb + \ee$.   

In the out-of-sample prediction problem for the linear model, the goal is to find rules for predicting unobserved future outcomes, $y_{new} = \x_{new}^T\bb + \e_{new}$, given the associated predictor vector $\x_{new}$ and the data $(\y,X)$.  In the formulation considered here, a prediction rule is determined by an estimator $\hat{\bb}$ for $\bb$ and the performance of the prediction rule is closely tied to properties of $\hat{\bb}$.  The ``usual'' estimator for $\bb$ is the ordinary least squares (OLS) estimator $\hat{\bb} = (X^TX)^{-1}X^T\y$.   However, the OLS estimator has drawbacks that are especially significant in high-dimensional data analysis, when the number of predictors $d$ is large, e.g. instability.  Furthermore, if $d >n$, then $X^TX$ is not invertible and the OLS estimator undefined (though this issue may be partially sidestepped by considering pseudoinverses, as is done below).  Thus, alternatives to the OLS estimator are desirable. 

Much of the recent research on high-dimensional linear models and alternatives to the OLS estimator has focused on sparsity.  
In this research, sparsity plays at least two roles: (i) sparse estimators for $\bb$ are often convenient, as they may aid interpretation\citep{tibshirani1996regression, fan2006statistical} and (ii) if $\bb$ is sparse in an appropriate sense, then this can often be leveraged to develop methods that perform very well, even with extremely high-dimensional datasets \citep{bunea2007sparsity, candes2007dantzig, bickel2009simultaneous, ye2010rate, zhang2010nearly, fan2011nonconcave, raskutti2011minimax, rigollet2011exponential}.  This provides a promising framework, which ideally yields interpretable estimators that perform well in high-dimensional data analysis.  However, several recent papers in genomics and statistics have questioned the degree of sparsity in modern genomic datasets (see, for instance, \citep{hall2009feature}, and the references contained therein -- including \citep{kraft2009genetic, goldstein2009common, hirschhorn2009genomewide} -- and, more recently, \citep{bansal2010statistical,manolio2010genomewide}).  This suggests that a closer study of non-sparse (or ``dense'') methods for high-dimensional linear models may prove useful.  

This paper contains a careful analysis of three non-sparse linear estimators that are alternative to the OLS estimator: ridge regression estimators, a class of James-Stein type estimators (scalar multiples of the OLS estimator), and marginal regression estimators.  We study the unconditional out-of-sample prediction error (predictive risk) associated with these estimators in a high-dimensional setting where the data are drawn from a multivariate normal distribution.  Though all of these estimators have been studied extensively in the past, the results in this paper offer unique and detailed insights into their comparative performance in high-dimensional data analysis, along with practical guidance for implementation and tuning parameter selection.  Symmetry properties of the estimators are easily leveraged in our formulation of the problem, which leads to many of the new insights delivered here.  No sparsity assumptions are made throughout the paper.  Though a direct comparative analysis of the estimators is emphasized, we also identify minimax ridge and James-Stein estimators (over the entire parameter space).  Broader optimality properties of non-sparse estimators are studied in \citep{dicker2012optimal}. 

Ultimately, the results in this paper have significant practical implications for high-dimensional linear models when little is known about the sparsity of the underlying signal, which may be partially summarized as follows: (i) if $\Cov(\x_i)$ is known or if a norm-consistent estimator is available, then the ridge estimator outperforms the James-Stein, OLS, and marginal estimators (in fact, results in \citep{dicker2012optimal} imply that the ridge estimator is nearly optimal for out-of-sample prediction in the described setting); (ii) if little is known about $\Cov(\x_i)$ or if $d/n$ is small, then the James-Stein estimator may be an effective alternative to the ridge and OLS estimators; and (iii) the marginal estimator has serious deficiencies for out-of-sample prediction.  

\section{Preliminaries: Definitions, notation, and an overview of results} 

\subsection{Out-of-sample prediction}

Each estimator, $\hat{\bb} = \hat{\bb}(\y,X)$, of $\bb$ determines a linear prediction rule, $\hat{y}(\x) = \x^T\hat{\bb}$.  We define the unconditional out-of-sample prediction error (predictive risk) of $\hat{\bb}$ to be
\begin{equation}\label{pr}
E\{y_{new} - \hat{y}(\x_{new})\}^2 = E(y_{new} - \x_{new}^T\hat{\bb})^2,
\end{equation}
where $(y_{new},\x_{new}^T)$ is independent of $(\y,X)$ and drawn from the same data-generating mechanism as $(y_i, \x_i^T)$, and the expectation in (\ref{pr}) is taken over $(y_{new},\x_{new})$ and $(\y,X)$.  The goal of the unconditional out-of-sample prediction problem is to minimize (\ref{pr}) over estimators $\hat{\bb}$.  

In order to evaluate (\ref{pr}), the distribution of $\e_i$ and $\x_i$ must be specified.  We assume that 
\begin{equation}\label{mvn1}
\x_1,...,\x_n \stackrel{\mbox{\tiny iid}}{\sim} N(0,\S) \mbox{ and } \e_1,...,\e_n \stackrel{\mbox{\tiny iid}}{\sim} N(0,\s^2) \mbox{ are independent,}
\end{equation}
where $\S$ is a $d \times d$ positive definite matrix and $\s^2 > 0$.  These distributional assumptions are restrictive.  However, other authors studying predictive risk have made similar assumptions \citep{stein1960multiple, baranchik1973inadmissibility, breiman1983many, brown1990ancillarity, leeb2009conditional} and we believe that the insights imparted by the resulting simplifications are worthwhile.

The assumptions (\ref{mvn1}) imply that $E(\x_i) = 0$ and $E(y_i) = 0$.  In other words, the model considered here does not have an intercept term.  In many practical settings, it is more appropriate to allow $E(y_i), E(\x_i) \neq 0$.  In fact, all of the methods studied in this paper can accomodate data with an intercept term, provided one first follows the usual approach of centering the data, and then decorrelating the observations (i.e., adjusting for the degree of freedom lost upon centering the data).  

Let $\w_i^T = (y_i,\x_i^T) \in \R^{d +1}$ and note that the assumption (\ref{mvn1}) is equivalent to assuming $\w_1,...,\w_n \stackrel{\mbox{\tiny iid}}{\sim}  N(0,V)$, where 
\begin{equation} \label{mvn2}
V = \left(\begin{array}{cc} \s^2 + \bb^T\S\bb & \bb^T\S \\ \S\bb & \S \end{array}\right) \in PD(d+1)
\end{equation}
and  $PD(d + 1)$ is the collection of all $(d+1)\times(d+1)$ positive definite matrices. 
The predictive risk (\ref{pr}) of an estimator $\hat{\bb}$ may be re-expressed as
\[
E(y_{new} - \x_{new}^T\hat{\bb})^2 = E_V\left\{(\hat{\bb} - \bb)^T\S(\hat{\bb} - \bb)\right\} + \s^2,
\]
where the subscript $V$ in the expectation on the right-hand side above indicates that the expectation is taken over $\w_1,...,\w_n \sim N(0,V)$.   After standardizing by $\s^2$, the predictive risk is equivalent to
\[
R_V(\hat{\bb}) = \s^{-2}E_V\left\{(\hat{\bb} - \bb)^T\S(\hat{\bb} - \bb)\right\}.
\]
In fact, $R_V(\hat{\bb})$ is the primary object of study in the sequel and we will typically refer to $R_V(\hat{\bb})$ itself as the predictive risk (or out-of-sample prediction error) of $\hat{\bb}$.  Note that the predictive risk $R_V(\hat{\bb})$ is completely determined by the estimator $\hat{\bb}$ and the positive definite matrix $V \in PD(d+1)$.  We will often write $E_{\S}(\cdot)$ in place of $E_V(\cdot)$ when the expectation only involves the random matrix $X$.  Similarly, we write $P_V(\cdot)$ or $P_{\S}(\cdot)$ when computing probabilities involving $\w_1,...,\w_n$ or $X$, respectively.  

\subsection{The estimators}

For a matrix $A$, let $A^-$ denote its Moore-Penrose pseudoinverse.  Below, we define the estimators studied in this paper.  All of the estimators, which have the form $\hat{\bb} = A\y$ for some $d \times n$ matrix $A$.   

\begin{center}
\begin{tabular}{ll}
{\em OLS estimator}: & $\hat{\bb}_{ols} = (X^TX)^{-}X^T\y$. \\
{\em James-Stein estimator}: & $\hat{\bb}_{js}(\l) = (1 + \l)^{-1}(X^TX)^{-}X^T\y, \ \ \l \geq 0$. \\
{\em Ridge regression estimator}: & $\hat{\bb}_r(\l) = (X^TX + n\l\S)^{-}X^T\y, \ \ \l \geq 0$. \\
{\em Marginal regression estimator}: & $\hat{\bb}_m = n^{-1}\S^{-}X^Ty$.
\end{tabular} 
\end{center}

\noindent {\em The OLS estimator}:  This version of the OLS estimator is defined for all $d,n$ because it utilizes the pseudoinverse $(X^TX)^-$.  

\noindent {\em The James-Stein estimator}:  A version of this estimator was proposed by \cite{stein1960multiple}.  The parameter $\l \geq 0$ is a tuning (or shrinkage) parameter that must be specified by the user.  \cite{baranchik1973inadmissibility} proved that for a certain data dependent $\hat{\l}_{bar}$, the estimator $\hat{\bb}_{bar} = \hat{\bb}_{js}(\hat{\l}_{bar})$ has smaller predictive risk than the OLS estimator (the estimator $\hat{\bb}_{bar}$ is discussed further in Section 8.2 below).  We refer to $\hat{\bb}_{js}(\l)$ as the James-Stein estimator because of its superficial resemblance to the James-Stein estimator for the normal means problem \citep{james1961estimation}.  Notice that $\hat{\bb}_{js}(\l)$ is a scalar multiple of the OLS estimator and $\hat{\bb}_{js}(0) = \hat{\bb}_{ols}$.  The James-Stein estimator is a shrinkage estimator and $\l$ determines the amount of shrinkage: $||\hat{\bb}_{\l}||$ is decreasing in $\l$, where $||\cdot||$ denotes the $\ell^2$-norm.       

\noindent {\em The ridge regression estimator}:  Many versions of the ridge estimator have been studied and have been shown to outperform the OLS estimator in a variety of settings \citep{tikhonov1943stability, hoerl1970ridge, golub1979generalized, casella1980minimax}.  Perhaps the most common version of the ridge estimator has the form $\hat{\bb}_0(\l) = (X^TX + n\l I)^{-1}X^T\y$.  The ridge estimator considered here $\hat{\bb}_r(\l)$ has convenient symmetry properties and can be derived from a class of generalized ridge estimators proposed by \cite{casella1980minimax} when considerations about out-of-sample prediction are taken into account.  Notice that $\hat{\bb}_r(\l)$ depends on the covariance matrix $\S = \Cov(\x_i)$.  In practice, if $\S$ is not known then it may be feasible to replace $\S$ with an estimate $\hat{\S}$ to obtain a modified ridge estimator.  The effect of replacing $\S$ with $\hat{\S}$ on prediction error is discussed in Section 8.1.   If prior information about $\S$ is available, then this can be incorporated into $\hat{\S}$, otherwise the sample covariance $\hat{\S} = n^{-1}X^TX$ may be used.  If $\hat{\S} = n^{-1}X^TX$ is used in place of $\S$, then the modified ridge estimator reduces to the James-Stein estimator, $\hat{\bb}_{js}(\l)$.  This was observed previously by \citet{oman1984different}.  Like the James-Stein estimator, the ridge estimator is a shrinkage estimator and  $\l \geq 0$ is a shrinkage parameter that must be specified by the user.
 
\noindent {\em Marginal regression estimator}: Variants of this estimator (that are often implemented with $\S = I$) are known to have desirable screening and variable selection properties and have been used extensively for related applications \citep{fan2008sure}.  Like the ridge estimator, the marginal estimator depends on the covariance matrix $\S$.  If $\S$ is not known, then it may be replaced with an estimate $\hat{\S}$.  Taking $\hat{\S} = n^{-1}X^TX$ gives the OLS estimator.  It seems reasonable to also consider linear shrinkage estimators based on the marginal estimator, such as $(1 + \l)^{-1}\hat{\bb}_m$, $\l \geq 0$; these estimators are discussed further in Section 8.3.

\section{Overview of results}

\subsection{Symmetry properties}

In addition to being linear, the estimators defined in the previous section have notable symmetry properties.  In particular, they are linearly equivariant and scale equivariant.  These properties help simplify predictive risk calculations and are discussed in Section 4.  More fundamentally, we argue in Section 4 that these are natural properties for non-sparse estimators.  

\subsection{Finite sample results}

Both finite sample and asymptotic properties are studied in this paper.   Finite sample results are the subject of Section 5.  In finite samples, we identify oracle ridge and James-Stein estimators, $\hat{\bb}_r^* = \hat{\bb}_r(\l_r^*)$ and $\hat{\bb}_{js}^* = \hat{\bb}_{js}(\l_{js}^*)$, that depend on the (typically unknown) signal-to-noise ratio
\begin{equation}\label{sn}
\eta^2 = \frac{\bb^T\S\bb}{\s^2}
\end{equation}
(Propositions 4-5).  These estimators have the smallest predictive risk among ridge and James-Stein estimators with non-random shrinkage parameters $\l \geq 0$.  Simplified formulas for the estimators' predictive risk are also obtained.  Similar results have been obtained in other settings, e.g. \citep{pinsker1980optimal}, \citep{goldenshluger2003optimal}.  The major novelty of our results is their simplicity and their applicability to out-of-sample prediction.  These results provide the means for an initial comparative analysis of the estimators (Section 5.5).  In particular, we show that   
\begin{equation}  \label{intro}
R_V(\hat{\bb}_r^*) \leq R_V(\hat{\bb}_{js}^*) < \left\{\begin{array}{c} R_V(\hat{\bb}_{ols}) \\ R_V(\hat{\bb}_m). \end{array}\right.
\end{equation}
To our knowledge, these are the first analytic results providing a direct comparison between the predictive risk of ridge regression and James-Stein estimators.  In Section 5.5 we also argue that the marginal estimator  $\hat{\bb}_m$ has serious deficiencies for out-of-sample prediction.

\subsection{Asymptotic results}

In Sections 6-7, we study asymptotic properties of the estimators in high-dimensional settings.  This helps provide a better understanding of the estimators' performance in high-dimensional settings.  Asymptotic regimes where $d/n \to 0$, $d/n \to \rho \in (0,\infty)$, and $d/n \to \infty$ are all considered.  

We find that in order to ensure consistency (i.e. asymptotically vanishing predictive risk) for any of the estimators considered here, one must have $d/n \to 0$.   This is a common feature of non-sparse estimators that can be framed more generally in terms of minimax problems over highly symmetric parameter spaces \citep{pinsker1980optimal,donoho1994minimax}.  Indeed, one way to formulate dense estimation and prediction problems is in terms of minimax problems over $\ell^2$-balls. \cite{dicker2012optimal} proved that the minimax rate for out-of-sample prediction over $\ell^2$-balls is proportional to $d/n$.  Thus, the estimators considered here achieve the minimax rate (this is not too noteworthy, as many estimators achieve the minimax rate for dense estimation and prediction; however, \cite{dicker2012optimal} also proved the stronger result that the ridge estimator $\hat{\bb}_r^*$ is asymptotically minimax over $\ell^2$-balls).  

Though the estimators studied here require $d/n \to 0$ for consistency, we devote much of our effort to studying asymptotic regimes where $d/n \to \rho > 0$.   Our interest in these regimes is motivated by the emergence of important problems in high-dimensional data analysis where the role of sparsity is unclear (like those cited in Section 1 above), which highlight the importance of characterizing and thoroughly understanding dense problems and estimators in settings where $d/n$ is significantly larger than 0.

\subsubsection{Related work: Sparse problems and ellipsoids}

In contrast with dense problems, in sparse problems it is known that consistent estimation and prediction may be possible even if $d/n \to \infty$ \citep{bickel2009simultaneous, bunea2007sparsity, candes2007dantzig, raskutti2011minimax, ye2010rate}.  However, the required sparsity conditions (e.g. $\ell^p$-sparsity, $0 \leq p < 2$ \citep{abramovich2006adapting}) may not hold in general and our motivating interest lies precisely in these situations.  Other conditions on $\bb$ may also allow for consistent estimation or prediction when $d$ is much larger than $n$; for instance, if $\bb$ belongs to an $\ell^2$-ellipsoid with decaying axes,  $B(c,\aa) = \{\bb \in \R^d; \ \a_1\b_1^2 + \cdots +  \a_d\b_d^2 \leq c^2\}$, where $c \geq 0$, $\aa = (\a_1,...,\a_d)^T \in \R^d$, and $0 \leq \a_1 \leq \cdots \leq \a_d$ \citep{pinsker1980optimal, goldenshluger2001adaptive, cavalier2002sharp, goldenshluger2003optimal}.  In this direction, Goldenshluger and Tsybakov's (2001, 2003) work may be most relevant to ours.  They study out-of-sample prediction with ``blockwise'' James-Stein estimators (which are, in a sense, a hybrid of the ridge and James-Stein estimators considered here) and obtain adaptive asymptotic minimax results over $\ell^2$-ellipsoids.  Ultimately, these estimators leverage asymmetry in ellipsoidal parameter spaces (e.g. rapidly increasing $\a_i$) to obtain faster rates of convergence.  In the highly symmetric case that is most relevant to our results, where the parameter space is an $\ell^2$-ball $B(c) = B(c,(1,...,1)^T)$, Goldenshluger and Tsybakov's results require $d/n \to 0$ to ensure consistency and do not apply if  $d/n \to \rho > 0$.   In general, ellipsoid conditions are natural for many inverse problems in nonparametric function estimation, but they may be overly restrictive in other settings, such as the genomic applications discussed in Section 1.   

\subsubsection{Oracle estimators}

Section 6 of this paper contains a detailed asymptotic analysis of the predictive risk of $\hat{\bb}_{ols}$, $\hat{\bb}_m$, $\hat{\bb}_{js}^*$, and $\hat{\bb}_r^*$.  We show that if $d/n \to 0$, then $R_V(\hat{\bb}_r^*)$ and $R_V(\hat{\bb}_{js}^*)$ are asymptotically equivalent; whether or not these estimators are asymptotically equivalent to the OLS estimator depends on magnitude of the signal-to-noise ratio (Proposition 9).  

The regime where $d/n \to \rho \in (0,\infty)$ appears to be the natural setting for studying the estimators considered in this paper.  Using results from random matrix theory, e.g. \citep{marcenko1967distribution}, in Section 6 we obtain closed-form expressions for the asymptotic predictive risk of $\hat{\bb}_r^*$, $\hat{\bb}_{js}^*$, $\hat{\bb}_{ols}$, $\hat{\bb}_m$ as $d/n \to \rho \in (0,\infty)$.  These formulas are new and imply that each of the estimators exhibits distinct behavior in this asymptotic regime.  In particular, the benefits of the ridge estimator over the James-Stein, OLS, and marginal estimators observed in finite samples (\ref{intro}) persist when we pass to the limit. This contrasts with the case where $d/n \to 0$ and the ridge and James-Stein estimators are asymptotically equivalent.  

Finally, if $d/n \to \infty$, then $\bb$ is non-estimable by any of the methods considered here, in the sense that their asymptotic performance is no better than that of $\hat{\bb}_{null} = 0$.  In fact, results in \citep{dicker2012optimal} imply that if $d/n \to \infty$, then $\hat{\bb}_{null}$ is in fact asymptotically minimax over $\ell^2$-balls.  

\subsubsection{Adaptive estimation}

Section 7 is concerned with adaptive estimation, when $d/n \to \rho \in (0,1)$.  The oracle estimators $\hat{\bb}_r^*$ and $\hat{\bb}_{js}^*$ depend on the signal-to-noise ratio $\eta^2$, which is typically unknown.  We show that if $d/n \to \rho \in (0,1)$, then $\eta^2$ may be replaced with an estimate $\hat{\eta}^2$ and that the resulting estimators (which are called adaptive estimators because they ``adapt'' to the signal-to-noise ratio) are asymptotically equivalent to the oracle estimators.  Note that this addresses the problem of tuning parameter selection for the ridge and James-Stein estimators.  A corollary of the main result in Section 7 (Corollary 2, Proposition 10) implies that the adaptive ridge and James-Stein estimators are minimax over $V \in PD(d+1)$, provided $0 < \th \leq d/n < \Th < 1$ for some constants $\th,\Th \in R$ and $d,n$ are sufficiently large.  The requirement $d < n$ for our results on adaptive estimation is related to the fact that if $d \geq n$, then $\y = X\hat{\bb}_{ols}$ and the usual estimator for $\s^2$, $\hat{\s}^2 = (n-d)^{-1}||\y - X\bb||^2$, is undefined.  It may be possible to utilize other estimators for $\s^2$ in settings where $d \geq n$, but this is not pursued in detail here.  

\subsection{Miscellanea}

Some additional topics are discussed in Section 8.  Recall that the ridge estimator $\hat{\bb}_r$ depends on the predictor covariance matrix $\Cov(\x_i) = \S$.  In Proposition 11, we show that if a norm-consistent estimator of the predictor covariance $\hat{\S}$ is available, then substituting $\hat{\S}$ for $\S$ does not affect the asymptotic predictive risk of the ridge estimator when $d \asymp n$.  We also show that a previously proposed minimax James-Stein estimator (Baranchik's estimator $\hat{\bb}_{bar}$, introduced in Section 2.2) is sub-optimal, in terms of asymptotic predictive risk when compared with the oracle and adaptive James-Stein estimators proposed in this paper (Proposition 12).  Finally, we consider a shrinkage estimator based on the marginal regression estimator and show that the ridge estimator outperforms this estimator in terms of predictive risk.   Section 9 contains a concluding discussion.  Unless explicitly stated otherwise, all propositions are proved in Appendix B.

\section{Linear equivariance and scale invariance}

In addition to being linear estimators, the estimators studied in this article share important symmetry properties.  

\vspace{.1in} 
{\em Definition 1.} An estimator $\hat{\bb} = \hat{\bb}(y,X,\S)$ is {\em linearly equivariant} if 
\begin{equation}\label{linequiv}
A^{-1}\hat{\bb}(\y,X,\S) = \hat{\bb}(\y,XA,A\S A^T)
\end{equation}
for all $d \times d$ invertible matrices $A$.  It is {\em scale invariant} if 
\[
\hat{\bb}(\y,X,\S) = \hat{\bb}(t\y,tX,t^2\S)
\]
for all scalars $t \in \R\setminus\{0\}$.  If an estimator is both linearly equivariant and scale invariant, we say that it is {\em LS}.  \hfill $\Box$

\vspace{.1in}

An estimator is linearly equivariant if it is compatible with linear transformations of the predictor basis; it is scale invariant if it is invariant under scaling of the data.  Note that estimators in Definition 1 are allowed to depend on the population predictor covariance $\S$ and the compatibility criterion (\ref{linequiv}) implies that a linearly equivariant estimator's dependence on $\S$ must respect changes of basis.  Speaking broadly, linearly equivariant estimators may be appropriate in situations where there is little prior knowledge about the information carried in the given predictor basis as it relates to the outcome.  By contrast, sparsity assumptions convey exactly this type of information and linear equivariance is less appropriate for sparse signals.  
 Notice that $\hat{\bb}_{ols}$ and $\hat{\bb}_m$, are LS; for any fixed $\l \geq 0$, $\hat{\bb}_{js}(\l)$ and $\hat{\bb}_r(\l)$ are also LS.  

The symmetry properties of the various classes of estimators described in Definition 1 lead to some useful simplifications in their predictive risk.  Recall the signal-to-noise ratio $\eta^2 = \bb^T\S\bb/\s^2$, defined in (\ref{sn}).

\vspace{.1in} 
{\em Proposition 1.} {\em(a)} If $\hat{\bb}$ is linearly equivariant, then $R_V(\hat{\bb}) = R_{V_0}(\hat{\bb})$, where
\[
V_0 = \left(\begin{array}{cc} \s^2 + \bb^T\S \bb & \bb^T\S^{1/2} \\ \S^{1/2}\bb & I \end{array}\right).
\]  

 {\em(b)} If $\hat{\bb}$ is LS, then $R_V(\hat{\bb}) = R_{V_{\u}}(\hat{\bb})$, where
\[
V_{\u} = \left(\begin{array}{cc} 1+ \eta^2 & \eta\u^T \\ \eta\u & I \end{array}\right),
\]
$\eta^2$ is the signal-to-noise ratio (\ref{sn}), and $\u \in \R^d$ is any fixed unit vector.  In particular, $R_V(\hat{\bb})$ depends on $\bb$, $\S$, and $\s^2$ only through the signal-to-noise ratio (and $d,n$).  \hfill $\Box$
\vspace{.1in} 

Proposition 1 implies that when computing the predictive risk of LS estimators, we may assume without loss of generality that $\S = I$, $\s^2 = 1$, and $\bb = \eta\u$, for an arbitrary fixed unit vector $\u \in \R^d$.  This is used repeatedly to prove the propositions below.  On the other hand, we do not make the blanket assumption that $\S = I$, $\s^2 = 1$, and $\bb = \eta^2\u$, in order to emphasize that some of the LS estimators considered here (namely, $\hat{\bb}_r(\l)$ and $\hat{\bb}_m$) depend on $\S$, while others ($\hat{\bb}_{ols}$ and $\hat{\bb}_{js}(\l)$) do not -- this distinction becomes less apparent if one assumes that $\S = I$.  

More fundamentally, Proposition 1 implies that for LS estimators and out-of-sample prediction,  sparsity is irrelevant.  Indeed, the risk of an LS estimator is completely determined by the signal-to-noise ratio $\eta^2 = \bb^T\S\bb/\s^2$, which does not capture well-accepted notions of sparsity (e.g. $\ell^p$-sparsity, $0 \leq p < 2$).  Thus, LS estimators are robust to sparsity assumptions, which may be desirable in situations where little is known about sparsity.  On the other hand, LS estimators are not able to take advantage of sparsity in situations where $\bb$ is in fact sparse.  

\section{Finite sample results}

\subsection{The OLS estimator} 

For $d  < n -1$, the predictive risk of the OLS estimator is well know.  For $d \geq  n-1$, the analysis is straightforward, but less widely available in the literature (recall that the Moore-Penrose pseudoinverse is used to define $\hat{\bb}_{ols}$ when $d > n$).  

\vspace{.1in} 
{\em Proposition 2.} 
\[
R_V(\hat{\bb}_{ols}) = \left\{\begin{array}{cl} \frac{d}{n-d-1} & \mbox{if } d < n-1 \\
\frac{n}{d - n - 1} +  \eta^2\frac{d-n}{d}& \mbox{if } d > n + 1 \\
\infty & \mbox{if } d \in \{n-1,n,n+1\}.  \end{array}\right.
\]\hfill $\Box$
\vspace{.1in}

Notice that the predictive risk of the OLS estimator is finite whenever $d \neq n-1,n,n+1$.  In particular, it is finite when $d > n + 1$; however, the OLS estimator is biased when $d > n$.  

\subsection{The marginal regression estimator}

To our knowledge, the predictive risk of the marginal regression estimator, $\hat{\bb}_m$ has not been studied previously. 

\vspace{.1in}

{\em Proposition 3.}  
\[
R_V(\hat{\bb}_m) = \frac{d}{n} + \eta^2\frac{d+1}{n}.
\] \hfill $\Box$

\subsection{The James-Stein estimator}

Proposition 4 yields the predictive risk of $\hat{\bb}_{\l}$, the optimal James-Stein tuning parameter $\l_{js}^*$, and the oracle James-Stein estimator $\hat{\bb}_{js}^* = \hat{\bb}_{js}(\l_{js}^*)$.  The result follows from  a bias-variance decomposition. 

\vspace{.1in}
{\em Proposition 4.} Let 
\[
\l_{js}^* = \left\{\begin{array}{cl} \frac{d}{\eta^2(n - d - 1)} & \mbox{if } d < n-1 \\
\frac{d}{\eta^2(d-n-1)} & \mbox{if } d > n + 1 \\
\infty & \mbox{if } d \in \{n-1,n,n+1\}\end{array}\right.
\]
and let $\hat{\bb}_{js}^* = \hat{\bb}_{js}(\l_{js}^*)$.  Then
\begin{equation}\label{prop3a}
R_V\{\hat{\bb}_{js}(\l)\} = \left\{\begin{array}{cl} \left(\frac{1}{1 + \l}\right)^2\frac{d}{n-d-1} + \left(\frac{\l}{1 + \l}\right)^2\eta^2 & \mbox{if } d < n-1 \\
\left(\frac{1}{1 + \l}\right)^2\frac{n}{d - n -1} + \left(\frac{\l}{1 + \l}\right)^2\frac{n}{d}\eta^2 + \frac{d-n}{d}\eta^2 &   \mbox{if } d > n + 1 \\ 
\infty & \! \! \begin{array}{l} \mbox{if } d \in \{n-1,n,n+1\} \\ \mbox{and } \l < \infty \end{array}. \end{array}\right.
\end{equation}
and
\[
\min_{\l \in [0,\infty]} R_V\{\hat{\bb}_{js}(\l)\} = R_V(\hat{\bb}_{js}^*) = \left\{ \begin{array}{cl} \frac{\eta^2d}{\eta^2(n - d - 1) + d} & \mbox{if } d < n-1 \\
\frac{\eta^2n}{\eta^2(d-n-1) + d} + \eta^2\frac{d-n}{d} & \mbox{if } d > n + 1 \\
\eta^2 & \mbox{if } d \in \{n-1, n,n+1\}. \end{array}\right.
\] \hfill $\Box$

\subsection{Ridge regression}

The predictive risk of $\hat{\bb}_r(\l)$, optimal ridge parameter $\l_r^*$, and the oracle ridge estimator are identified in Proposition 5.  We have been unable to find a closed form expression for the predictive risk of the oracle ridge estimator -- such an expression may not exist.  

\vspace{.1in} 

{\em Proposition 5. } Let $\l_r^* = d/(n\eta^2)$ and let $\hat{\bb}_r^* = \hat{\bb}_r(\l_r^*)$.  Then
\begin{equation}\label{prop6a}
R_V\{\hat{\bb}_r(\l)\} = E_I\tr \left\{(X^TX + n\l I)^{-2}\left(X^TX + \frac{\eta^2n^2}{d}\l^2 I\right)\right\}.
\end{equation}
and
\begin{equation}\label{prop6b}
R_V(\hat{\bb}_r^*) = \inf_{\l \in [0,\infty]} R_V\{\hat{\bb}_r(\l)\} = E_I\tr(X^TX + n\l_r^*I)^{-1}.
\end{equation} \hfill $\Box$

 \subsection{A comparative analysis, part I: Oracle estimators and finite sample predictive risk}

Propositions 2, 4, and 5 immediately imply that $R_V(\hat{\bb}_{js}^*)$, $R_V(\hat{\bb}_r^*) < R_V(\hat{\bb}_{ols})$.  The predictive risk of $\hat{\bb}_{ols}$, $\hat{\bb}_m$, $\hat{\bb}_{js}^*$, and $\hat{\bb}_r^*$ are all increasing in $\eta^2$.  Furthermore, if $d \notin \{n - 1,n,n+1\}$, then
\begin{equation}\label{ratlim}
\lim_{\eta^2 \to \infty} \frac{R_V(\hat{\bb}_{js}^*)}{R_V(\hat{\bb}_{ols})} = \lim_{\eta^2 \to \infty} \frac{R_V(\hat{\bb}_r^*)}{R_V(\hat{\bb}_{ols})} = 1.
\end{equation}
On the other hand, if $\eta^2 = 0$, then
\[
R_V(\hat{\bb}_{js}^*) = R_V(\hat{\bb}_r^*) = 0 < R_V(\hat{\bb}_{ols}).
\]
These observations may be summarized as follows: (i) shrinkage estimators offer improvements over the OLS estimator in terms of out-of-sample prediction and (ii) these improvements are most substantial when the signal-to-noise ratio $\eta^2$ is small and diminish as $\eta^2$ grows larger.  Properties like these are common among shrinkage estimators in many contexts.  

At this stage, it appears to be difficult to make a detailed comparison between the predictive risk of the James-Stein estimator and that of the ridge estimator.  However, we have the following result.

\vspace{.1in}
{\em Proposition 6.} For any $\bb \in \R^d$, 
\[
R_V(\hat{\bb}_r^*) \leq R_V(\hat{\bb}_{js}^*) < R_V(\hat{\bb}_{ols}).
\]
The inequality on the left is strict unless $\eta^2=0$.  \hfill $\Box$

Proposition 6 is a consequence of Jensen's inequality and implies that the oracle ridge estimator $\hat{\bb}_r^*$ has smaller predictive risk than the oracle James-Stein estimator $\hat{\bb}_{js}^*$.  This is, perhaps, not surprising, given that the ridge estimator $\hat{\bb}_r^*$ utilizes knowledge of the predictor covariance $\S$ while the James-Stein estimator $\hat{\bb}_{ols}^*$ does not.  On the other hand, the marginal estimator $\hat{\bb}_m$ also utilizes knowledge of the predictor covariance, but the next proposition suggests that it is less suitable for out-of-sample prediction.   

\vspace{.1in}

{\em Proposition 7.} Let $\hat{\bb}_{null} = 0$. 
\begin{itemize}
\item[(a)] Suppose that $d < n - 1$.  Then 
\begin{eqnarray*}
R_V(\hat{\bb}_m) \leq R_V(\hat{\bb}_{ols}) &\mbox{ if and only if } &\eta^2 \leq \frac{d}{n-d-1} \\
& \mbox{if and only if} & R_V(\hat{\bb}_{null}) \leq R_V(\hat{\bb}_m).
\end{eqnarray*}
\item[(b)] Suppose that $d \geq n - 1$.  Then $R_V(\hat{\bb}_{null}) \leq R_V(\hat{\bb}_m)$. 
\end{itemize}
\hfill $\Box$

\vspace{.1in}

The proof of Proposition 7 is straightforward and is omitted.  Proposition 7 implies that if the marginal estimator has smaller predictive risk than the OLS estimator, then
the marginal estimator itself is outperformed by the trivial estimator $\hat{\bb}_{null} = 0$.  Additional drawbacks of the marginal estimator include that for any fixed $d,n$, $\lim_{\eta^2 \to \infty} R_V(\hat{\bb}_m) = \infty$.  On the other hand, if $d < n- 1$, $\lim_{\eta^2 \to \infty} R_V(\hat{\bb}_{ols}) = \lim_{\eta^2 \to \infty} R_V(\hat{\bb}_{js}) = \lim_{\eta^2 \to \infty} R_V(\hat{\bb}_{r}) = d/(n-d-1)$ (if $d > n$, then the limiting risk of these estimators is infinite as well; this is related to non-estimability issues that are discussed further in Section 6.2.3).  

A direct corollary of Proposition 7 is that $\hat{\bb}_m$ is dominated by the estimator 
\begin{equation}\label{dom}
\hat{\bb}_{dom} = \left\{\begin{array}{cl} 0 & \mbox{if } \eta^2 < \frac{d}{n-d-1} \mbox{ or } d \geq n-1 \\ \hat{\bb}_{ols} &\mbox{otherwise} \end{array}\right.  
\end{equation}
in the sense that $R_V(\hat{\bb}_{dom}) \leq R_V(\hat{\bb}_m)$ with strict inequality whenever $\eta^2 \neq d/(n-d-1)$.  Since 
\[
R_V(\hat{\bb}_{js}^*) \leq \eta^2 = R_V(\hat{\bb}_{null}),
\]
 with strict inequality unless $\eta^2 = 0$, it follows that $R_V(\hat{\bb}_{js}^*) \leq R(\hat{\bb}_{dom})$ with equality if and only if $\eta^2 = 0$.  Thus, we recover (\ref{intro}):
\[
R_V(\hat{\bb}_r^*) \leq R_V(\hat{\bb}_{js}^*) < \left\{\begin{array}{l} R_V(\hat{\bb}_{ols}) \\ R_V(\hat{\bb}_m). \end{array}\right.
\]

\section{Asymptotic results for the oracle estimators}

The rest of this paper is primarily concerned with asymptotic results, which can be divided into two categories:  results about the asymptotic predictive risk of oracle estimators and results about how well adaptive estimators approximate the oracle estimators.  Asymptotic properties of oracle estimators are studied in this section.  

\subsection{Ridge regression}

To obtain a formula for the asymptotic predictive risk of $\hat{\bb}_r^*$, we rely on classical results from random matrix theory.  For $\rho \in (0,\infty)$ the Mar\v{c}enko-Pastur density $f_{\rho}$ is defined by
\[
f_{\rho}(z)  = \frac{dF_{\rho}}{dz}(z)  = \max\left\{1 - \rho^{-1}, 0\right\}\d_{0}(z) + \frac{1}{2\pi\rho z}\sqrt{4\rho - (z - \rho - 1)^2}, \ \ \ a \leq z \leq b,
\]
where $a =  (1 - \sqrt{\rho})^2$, $b = (1 + \sqrt{\rho})^2$, and $\d_0(x) = 1$ or $0$ according to whether $x = 0$ or $x \neq 0$.  The density $f_{\rho}$ determines the Mar\v{c}enko-Pastur distribution, $F_{\rho}$, which is the limiting distribution of the eigenvalues of $n^{-1}X^TX$, if $\S = I$, $n \to \infty$, and $d/n \to \rho \in (0,\infty)$ \citep{marcenko1967distribution}.  The Stieltjes transform of the Mar\v{c}enko-Pastur distribution,
\begin{equation}  
m_{\rho}(s)  =  \int \frac{1}{z - s} \ dF_{\rho}(z) 
 =  -\frac{1}{2\rho s}\left\{s + \rho - 1 +  \sqrt{(s + \rho - 1)^2 - 4\rho s}\right\}, \ \ s < 0, \label{stieltjes}
\end{equation}
\citep{bai1993convergence} has played a prominent role in the discovery and subsequent analysis of the Mar\v{c}enko-Pastur distribution; see, for instance, \citep{bai1993convergence}, \citep{silverstein1995strong}, and \citep{elkaroui2008spectrum}. The main result of this section implies that the risk of the oracle ridge estimator $R_V(\hat{\b}_r^*)$ may be approximated by $(d/n)m_{d/n}(-\l_r^*/n)$, where $\l_r^*$ is the optimal ridge parameter defined in Proposition 5.

\vspace{.1in}
{\em Proposition 8.}
Suppose that $0 < \th \leq d/n \leq \Th < \infty$ for some fixed constants $\th,\Th \in \R$. 
\begin{itemize}
\item[(a)] If $0 < \th < \Th < 1$ or $1 < \th < \Th < \infty$ and $n - d > 5$, then
\[
\left| R_V(\hat{\bb}_r^*) - \frac{d}{n}m_{d/n}(-\l_r^*)\right| = O\left(\frac{\eta^2}{1 + \eta^2}n^{-1/4}\right).
\]
\item[(b)] If $0 < \th < 1 < \Th < \infty$, then
\[
\left| R_V(\hat{\bb}_r^*) - \frac{d}{n}m_{d/n}(-\l_r^*)\right| = O(\eta^2 n^{-5/48}).
\]
\end{itemize} \hfill $\Box$
\vspace{.1in}

There are two keys to the proof of Proposition 8. The first is the observation that 
\[
\frac{n}{d}R_V(\hat{\bb}^*_r) = \frac{1}{d} E_I\tr(n^{-1}X^TX + \l_r^*I)^{-1} = E_I \int \frac{1}{s + \l_r^*} \ d\F_{n,d}(s),
\]
where $\F_{n,d}$ is the empirical cumulative distribution function of the eigenvalues of $n^{-1}X^TX$ -- in other words, the risk of the oracle ridge estimator is the expected value of the Stieltjes transform of $\F_{n,d}$.  The second key is Theorem 3.1 of \citet{bai1993convergence} which states that
\begin{equation}\label{bai}
\sup_s |E\F_{n,d}(s) - F_{d/n}(s)| = \left\{\begin{array}{cl} O(n^{-1/4}) &  \mbox{if } 0 < \th < \Th < 1 \mbox{ or}  1 < \th < \Th < \infty, \\ O(n^{-5/48}) & \mbox{if } 0 < \th < 1 < \Th < \infty. \end{array}\right. 
\end{equation}
The different rates in (\ref{bai}) for settings where $1 < \th < \Th < \infty$ or $1 < \th < \Th < \infty$ and $0 < \th < 1 < \Th < \infty$ helps to explain why these situations are considered separately in Proposition 8.  

\vspace{.1in}
{\em Corollary 1.}
Define the asymptotic predictive risk of the oracle ridge estimator, 
\[
R_r(\rho,\eta^2)  =  \rho m_{\rho}(-\rho/\eta^2)   =  \frac{1}{2\rho}\left[\eta^2(\rho - 1) - \rho + \sqrt{\{\eta^2(\rho - 1) - \rho\}^2 + 4\rho^2\eta^2}\right].
\]
Then
\begin{equation}\label{cor1a}
\lim_{d/n \to \rho} \sup_{V \in PD(d+1)} \left|R_V(\hat{\bb}_r^*) - R_r(d/n,\eta^2)\right| = 0,
\end{equation}
provided $\rho \in (0,\infty)\setminus\{1\}$.  \hfill $\Box$  
\vspace{.1in}

{\em Remark.} The limit (\ref{cor1a}) indicates that $n \to \infty$ and $d/n \to \rho$.  

\subsection{A comparative analysis, part II: Asymptotic predictive risk}

\subsubsection{$d/n \to 0$}

Propositions 2-5 imply $R_V(\hat{\bb}_{ols})$, $R_V(\hat{\bb}_{js}^*)$, $R_V(\hat{\bb}_{js}^*) = O(\eta^2d/n)$ and $R_V(\hat{\bb}_m^*) = O\{(1 + \eta^2)d/n\}$.  It follows that if $d/n \to 0$, then the estimators are consistent.  Additionally, we have the following result.

\vspace{.1in}
{\em Proposition 9.}
If $d/n \to 0$ and $d/(n\eta^2) \to c \in [0,\infty]$, then
\[
\frac{R_V(\hat{\bb}_{js}^*)}{R_V(\hat{\bb}_{ols})} \to \frac{1}{1 + c} \mbox{ and } \frac{R_V(\hat{\bb}_{r}^*)}{R_V(\hat{\bb}_{js}^*)} \to 1.
\] \hfill $\Box$
\vspace{.1in}

Thus, if $d/n \to 0$ and the signal-to-noise ratio is large (i.e. $d/(n\eta^2) \to0$), then the OLS, James-Stein, and ridge estimators are all asymptotically equivalent; if $d/n \to 0$ and the signal-to-noise ratio is small ($d/(n\eta^2) \to c > 0$), then the James-Stein and ridge estimators are asymptotically equivalent and they both outperform the OLS estimator asymptotically.  

\subsubsection{$d/n \to \rho \in (0,\infty)$}

If $d/n \to \rho \in (0,\infty)$, then Corollary 1 implies that the predictive risk of $\hat{\bb}_r^*$ is non-vanishing.  This is also true of the OLS, James-Stein, and marginal estimators.  In fact, it is straightforward to derive the asymptotic predictive risk of these estimators (all of the limits below are valid for fixed signal-to-noise ratios $\eta^2$; in fact, the convergence holds for varying degrees of uniformity in $\eta^2$ for the different estimators, however, these details are not critical for our analysis).
\begin{eqnarray*}
\mbox{\em OLS}: \ \  R_{ols}(\rho,\eta^2)  & = & \lim_{d/n \to \rho} R_V(\hat{\bb}_{ols}) \ \  =  \ \     \left\{\begin{array}{cl} \frac{\rho}{1 - \rho} & \mbox{if } \rho < 1 \\ \frac{1}{\rho - 1} + \eta^2\frac{\rho - 1}{\rho} & \mbox{if } \rho > 1 \\ \infty & \mbox{if } \rho = 1,\end{array}\right. \\
\mbox{\em James-Stein}: \ \  R_{js}(\rho,\eta^2) &  = &   \lim_{d/n \to \rho} R_V(\hat{\bb}_{js}) \ \ = \ \    \frac{\eta^2(\rho \wedge 1)}{\eta^2|1 - \rho| + \rho} + \eta^2\frac{(\rho \vee 1) - 1}{\rho} \\
\mbox{\em Marginal}: \ \  R_m(\rho,\eta^2) &  = &  \lim_{d/n \to \rho} R_V(\hat{\bb}_{m}) \ \ = \ \    (\eta^2 + 1)\rho.
\end{eqnarray*}
It is easy to check that
\begin{equation} \label{asympineq} 
R_r(\rho,\eta^2) \leq R_{js}(\rho,\eta^2) < \left\{\begin{array}{l} R_{ols}(\rho,\eta^2) \\ R_m(\rho,\eta^2), \end{array}\right.
\end{equation}
and that the inequality on the left is strict unless $\eta^2 = 0$.  The inequalities (\ref{asympineq}) indicate that the advantages of the ridge estimator (over the James-Stein, marginal, and OLS estimators) and the James-Stein estimator (over the marginal and OLS estimators) persist in high-dimensional datasets, as $n \to \infty$ and $d/n \to \rho \in (0,\infty)$.  More fundamentally, they illustrate that different linear estimators may have significantly different out-of-sample prediction properties in high-dimensional data analysis -- differences between the estimators do not ``wash-out'' in this asymptotic setting.  In addition to providing this qualitative information, these asymptotic formulas provide an analytic tool for studying the various estimators' predictive risk.  Figures 1-3 contain several plots of asymptotic predictive risk for the OLS, marginal, oracle James-Stein, and oracle ridge estimators.   We point out that Figures 2-3 contain plots of the asymptotic predictive risk for the oracle marginal shrinkage estimator, which is introduced in Section 8.3 (Proposition 13).  

\begin{figure}[h]
\begin{center}
\includegraphics[width=\textwidth]{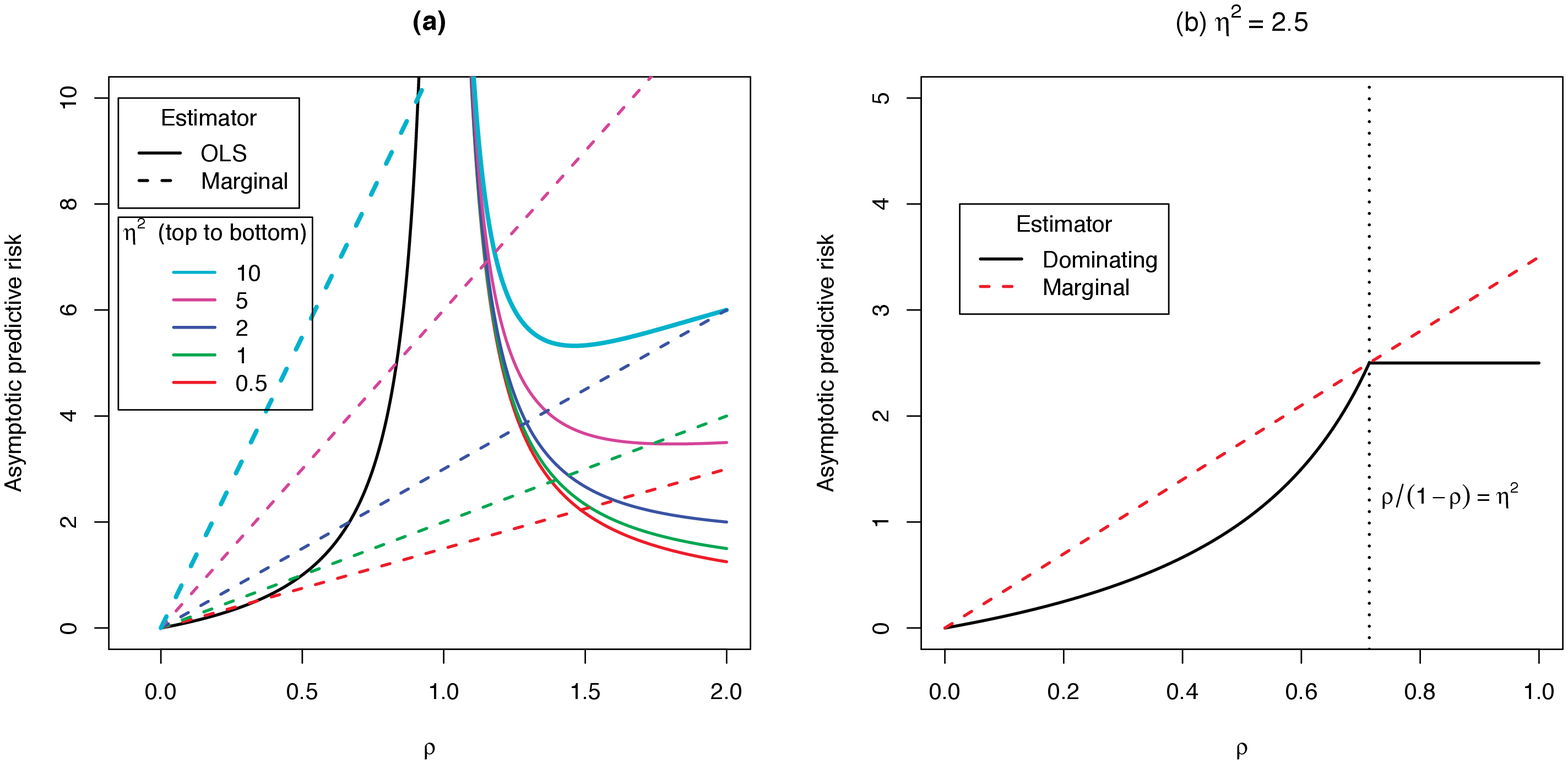}
\end{center}
\caption{ (a) The asymptotic predictive risk of the OLS ($\hat{\bb}_{ols}$) and marginal ($\hat{\bb}_m$) estimator versus $\rho$ for multiple values of $\eta^2$. The asymptotic predictive risk is increasing with $\eta^2$.  (b) The asymptotic predictive risk of the marginal estimator and the dominating estimator ($\hat{\bb}_{dom}$, defined in (\ref{dom})) versus $\rho$, for $\eta^2 = 2.5$. One easily checks that $\lim_{d/n \to \rho} R_V(\hat{\bb}_{dom}) = \rho/(1 - \rho)$ or $\eta^2$, according to whether $\rho/(1 - \rho) < \eta^2$ or $\rho/(1 - \rho) \geq \eta^2$.  }
\end{figure}
\begin{figure}[h]
\begin{center}
\includegraphics[width=\textwidth]{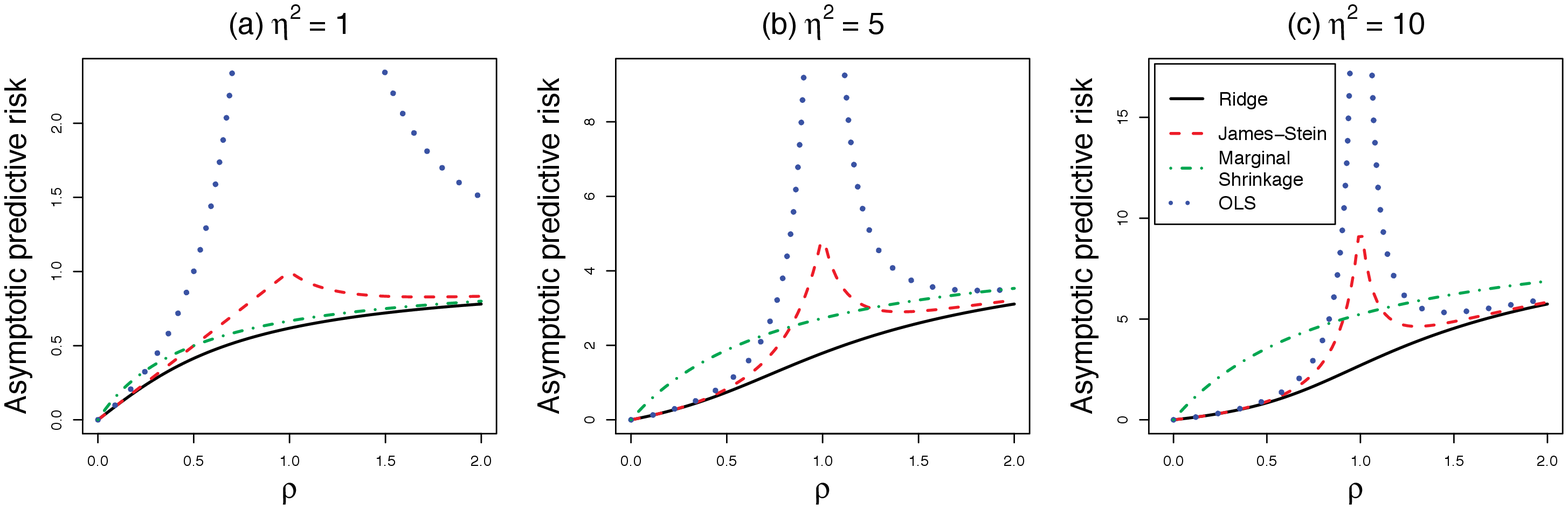}
\end{center}
\caption{Asymptotic predictive risk versus $\rho$ for the oracle ridge estimator ($\hat{\bb}_r^*$), the oracle James-Stein estimator ($\hat{\bb}_{js}^*$), the oracle marginal shrinkage estimator ($\hat{\bb}_m^*$, defined in Proposition 13), and the OLS ($\hat{\bb}_{ols}$) estimator for various values of $\eta^2$: (a) $\eta^2 = 1$, (b) $\eta^2 = 5$, and (c) $\eta^2 = 10$.  }
\end{figure}
\begin{figure}[h]
\begin{center}
\includegraphics[width=\textwidth]{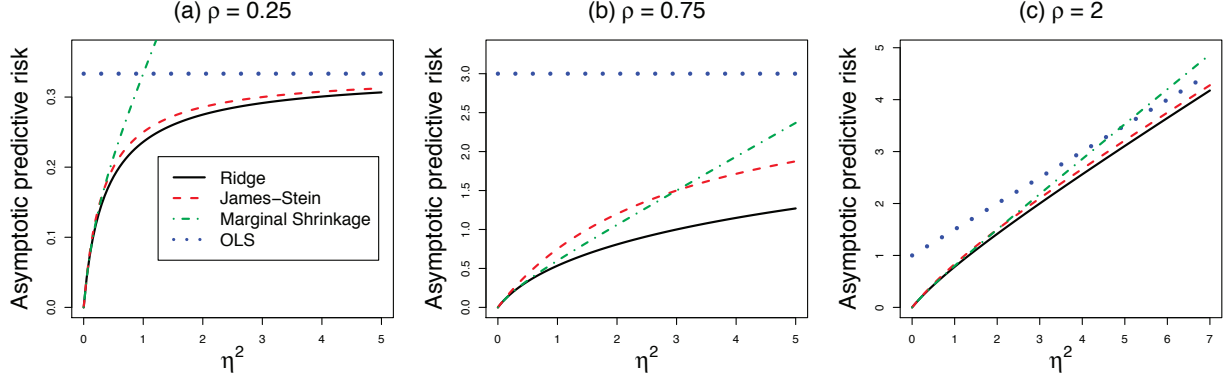}
\end{center}
\caption{Asymptotic predictive risk versus the signal-to-noise ratio, $\eta^2$,   for the oracle ridge estimator ($\hat{\bb}_r^*$), the oracle James-Stein estimator ($\hat{\bb}_{js}^*$),  the oracle marginal shrinkage estimator ($\hat{\bb}_m^*$, defined in Proposition 13), and the OLS ($\hat{\bb}_{ols}$) estimator for various values of $\rho$: (a) $\rho = 0.25$, (b) $\rho = 0.75$, and (c) $\rho = 2$.}
\end{figure}

Figure 2 depicts the singularities in $R_{ols}(\rho,\eta^2)$ and $R_{js}(\rho,\eta^2)$ at $\rho = 1$.   Notice that the ridge estimator's advantage over the other estimators appears to be most pronounced at $\rho = 1$.  This is borne out by the fact that as $\eta^2 \to \infty$, $R_r(1,\eta^2) \asymp \eta$, but $R_{js}(1,\eta^2) = \eta^2$.  On the other hand, if $\rho \neq 1$, then  $ R_r(\rho,\eta^2)/R_{js}(\rho,\eta^2) \to 1$ as $\eta^2 \to \infty$. 

\subsubsection{$d/n \to \infty$}

Figure 2 also suggests that for fixed $\eta^2$ each of the depicted estimator's asymptotic predictive risk approaches the same finite limit, as $\rho \to \infty$.   One can easily check that this limit is $\eta^2$.  This is reflective of the estimators' behavior as $d/n \to \infty$.  Indeed, it follows directly from Propositions 2-5 that 
\[
\lim_{d/n \to \infty} R_V(\hat{\bb}_{ols}) = \lim_{d/n \to \infty} R_V(\hat{\bb}_{js}^*) = \lim_{d/n \to \infty} R_V(\hat{\bb}_r^*) = \lim_{d/n \to \infty} R_V(\hat{\b}_m^*) = \eta^2
\]
($\hat{\bb}_m^*$ is the marginal shrinkage estimator defined in Proposition 13; for the marginal estimator, $R_V(\hat{\bb}_m) \to \infty$ as $\rho \to \infty$).  In fact, it can be shown that if $d > n$, then the predictive risk of any LS estimator must be at least $\eta^2(d-n)/d$.  Thus, if $\hat{\bb}$ is an LS estimator, then $\liminf_{d/n \to \infty} R_V(\hat{\bb})/R_V(\hat{\bb}_{null}) \geq 1$.  This is discussed further in \citep{dicker2012optimal}, where it is shown that if $d/n \to \infty$, then $\hat{\bb}_{null}$ is asymptotically minimax over $\ell^2$-balls.

\section{Adaptive estimators}

In the previous sections, we studied the predictive risk of $\hat{\bb}_r^*$ and $\hat{\bb}_{js}^*$. This analysis provides substantial insight into the performance of ridge and James-Stein estimators.   However, even assuming that $\S$ is known, the estimators $\hat{\bb}_r^*$ and $\hat{\bb}_{js}^*$ are usually not implementable, since they depend on the signal-to-noise ratio $\eta^2 = \bb^T\S\bb/\s^2$, which is usually unknown.  In this section, we show that if $0 < \inf d/n \leq \sup d/n < 1$, then $\eta^2$ may be effectively estimated.   More specifically, we show that the predictive risk of adaptive ridge and James-Stein estimators that utilize an estimate of $\eta^2$ (and, subsequently, an estimate of the optimal shrinkage parameter) is very close to that of the corresponding oracle estimators.  

\subsection{Estimating the signal-to-noise ratio}

For $d < n$, define the estimator
\begin{equation}\label{esnr}
\hat{\eta}^2 = \max\left\{\frac{||\y||^2}{n\hat{\sigma}^2} - 1, 0\right\} = \max\left\{\frac{||X\hat{\bb}_{ols}||^2}{n\hat{\sigma}^2} - \frac{d}{n}, 0\right\},
\end{equation}
where $\hat{\sigma}^2 = (n-d)^{-1}||\y - X\hat{\bb}_{ols}||^2$.  To motivate this definition, notice that
\[
\frac{||\y||^2}{n\hat{\sigma}^2} - 1 \approx \frac{\bb^T\S\bb + \s^2}{\s^2} - 1 = \eta^2.  
\] 
Results in Appendix A establish convergence rates for $E_V|\hat{\eta}^2 - \eta^2|^k$ and other technical results that are important for proving Proposition 10 below.  When $d \geq n$, $\hat{\s}^2$ is undefined and the estimator $\hat{\eta}^2$ breaks down.  We conjecture that it is possible to derive effective estimators of $\eta^2$ when $d \geq n$, provided $\sup d/n < \infty$.  This is discussed further in Section 9.1, however, it is not pursued at length in this paper.  

\subsection{The adaptive ridge and James-Stein estimators}

Recall from Propositions 4 and 5 that the oracle ridge and James-Stein estimators are $\hat{\bb}_r^* = \hat{\bb}_r(\l_r^*)$ and $\hat{\bb}_{js}^* = \hat{\bb}_{js}(\l_{js}^*)$, respectively, where $\l_r^* = d/(n\eta^2)$ and $\l_{js}^* = d/\{\eta^2(n - d - 1)\}$.  The next proposition is the main result in this section. 

\vspace{.1in}
{\em Proposition 10.}
Suppose that $0 < \th < d/n \leq \Th < 1$ for some fixed constants $\th, \Th \in\R$ and that $n - d > 6$.
\begin{itemize}
\item[(a)] Define the adaptive ridge estimator $\check{\bb}_r = \hat{\bb}_r(\hat{\l}_r^*)$, where $\hat{\l}_r^* = d/(n\hat{\eta}^2)$.  Then 
\[
R_V(\check{\bb}_r) = R_V(\hat{\bb}_r^*) +O\left\{\frac{1}{\sqrt{n}(\eta^2 + 1)}\right\}.
\]  
\item[(b)] Define the adaptive James-Stein estimator $\check{\bb}_{js} = \hat{\bb}_{js}(\hat{\l}_{js}^*)$, where $\hat{\l}_{js}^* = d/\{\hat{\eta}^2(n - d - 1)\}$.  Then
\[
R_V(\check{\bb}_{js}) = R_V(\hat{\bb}_{js}^*) +O\left\{\frac{1}{\sqrt{n}(\eta^2 + 1)}\right\}.
\]\hfill $\Box$
\end{itemize} 
\vspace{.1in}

{\em Remark 1.} When $\hat{\eta}^2 = 0$, we follow that convention that $\check{\bb}_r = \check{\bb}_{js} = 0$.  Notice that both $\check{\bb}_r$ and $\check{\bb}_{js}$ are LS.  

{\em Remark 2.} Proposition 10 implies that the predictive risk of the adaptive ridge and James-Stein estimators converge to the predictive of the corresponding oracle estimators uniformly for $V \in PD(d+1)$.  Moreover, if $\eta^2 \gg n^{-1/2}$, then Proposition 10 implies that 
\[
\frac{R_V(\check{\bb}_r)}{R_V(\hat{\bb}_{r}^*)}, \ \frac{R_V(\check{\bb}_{js})}{R_V(\hat{\bb}_{js}^*)} \to 1.  
\]
On the other hand, if $\eta^2 = O(n^{-1/2})$, then Proposition 10 is less useful for comparing the performance of $\check{\bb}_r$ and $\check{\bb}_{js}$ to the oracle estimators.  Indeed,  if $\eta^2 = O(n^{-1/2})$, then $R_V(\check{\bb}_{js})$, $R_V(\check{\bb}_{r})$, $R_V(\hat{\bb}_{js}^*)$, $R_V(\hat{\bb}_{r}^*)$, $R_V(\hat{\bb}_{null}) = O(n^{-1/2})$ (recall that $\hat{\bb}_{null} = 0$) and the benefits of $\check{\bb}_r$ and $\check{\bb}_{js}$ over even $\hat{\bb}_{null}$ are unclear; however, in this setting, the estimators are consistent and $R_V(\hat{\bb}_{null})/R_V(\hat{\bb}_r^*)$, $R_V(\hat{\bb}_{null})/R_V(\hat{\bb}_r^*) = O(1)$, so that even the oracle ridge and James-Stein estimators are not dramatic improvements on $\hat{\bb}_{null}$.

{\em Remark 3.} The condition $0 < \th \leq d/n$ in Proposition 10 can be removed.  However, the corresponding error terms in part (a) and (b) are more complicated; a precise statement is omitted from this paper.  Ultimately, however, when $d/n \to 0$ the message remains the same: if $\eta^2$ is not too small, then the adaptive estimators perform nearly as well as the oracle estimators. 

Proposition 10 is proved in Appendix B.  

\subsection{A comparative analysis, part III: Minimax estimators}

In addition to providing direct information about the performance the adaptive James-Stein and ridge estimators, vis-\`a-vis the oracle estimators, Proposition 10 also helps to shed light on the performance of the adaptive estimators relative to the OLS estimator and to each other.  Consider the following corollary to Proposition 10.  

\vspace{.1in}
{\em Corollary 2.}
Suppose that $0 < \th \leq d/n \leq \Th < 1$ for some fixed constants $\th, \Th \in \R$ and let $\check{\bb}_r$ and $\check{\bb}_{js}$ be the adaptive ridge and James-Stein estimators defined in the statement of Proposition 10.  If $n$ is sufficiently large, then
\begin{equation}\label{cor2a}
R_V(\check{\bb}_r), R_V(\check{\bb}_{js}) < R_V(\hat{\bb}_{ols}) \mbox{ for all } V \in PD(d + 1).
\end{equation}
In particular, if $n$ is sufficiently large, then  the adaptive ridge and James-Stein estimators are minimax over the entire parameter space in the sense that
\begin{equation}\label{cor2b}
\sup_{V \in PD(d+1)} R_V(\check{\bb}_r) = \sup_{V \in PD(d+1)} R_V(\check{\bb}_{js}) = \inf_{\hat{\bb}} \sup_{V \in PD(d + 1)} R_V(\hat{\bb}),
\end{equation}
where the infimum on the right-hand side of (\ref{cor2b}) is taken over all measurable estimators $\hat{\bb}$.    \hfill $\Box$
\vspace{.1in}

The first part of Corollary 2 (regarding the inequalities (\ref{cor2a})) follows directly from Proposition 10 and two observations: (i) if $0 < \th \leq d/n \leq \Th < 1$ for constants $\th,\Th \in \R$, then there exists a constant $c >0$ such that
\[
R_V(\hat{\bb}_{ols}) - R_V(\hat{\bb}_{js}^*) > \frac{c}{\eta^2},
\]
whenever $\eta^2$ and $n$ are sufficiently large, and (ii) $R_V(\hat{\bb}_r^*) \leq R_V(\hat{\bb}_{js}^*)$.  The minimaxity result (\ref{cor2b}) follows directly from the first part of the corollary and standard arguments from decision theory.

\section{Additional topics}

\subsection{The ridge estimator:  Estimating the population covariance}

If $\S$ is unknown, then even the adaptive ridge estimator $\check{\bb}_r$ can not be implemented.  On the other hand, if $\hat{\S}$ is an estimator of $\S$, then it may be reasonable to use a ridge estimator of the form
\[
\tilde{\bb}_r = \tilde{\bb}_r(\hat{\l}_r^*,\hat{\S}) = (X^TX + n\hat{\l}_r^*\hat{\S})^{-1}X^T\y,
\]
where $\hat{\l}_r^*$ is defined in Proposition 10.  For a matrix $A$, let $||A||$ denote it operator norm.  

\vspace{.1in}
{\em Proposition 11.}
Suppose that $0 < \th \leq d/n \leq \Th < 1$ for some fixed constants $\th, \Th \in \R$ and that $n - d > 6$.  Then
\[
R_V(\tilde{\bb}_r) = R_V(\hat{\bb}_r^*) + O\left\{||\S^{-1}||\ \left(E_V||\hat{\S} - \S||^2\right)^{1/2}\right\}  + O\left\{\frac{1}{(\eta^2+1)n^{1/2}}\right\}.
\]\hfill $\Box$
\vspace{.1in}

Proposition 11 implies that if the smallest eigenvalue of $\S$ is bounded below by some positive number and if $\hat{\S}$ is operator norm-consistent in the sense that $E_V||\hat{\S} - \S||^2 \to 0$, then the asymptotic predictive risk of $\tilde{\bb}_r = \tilde{\bb}_r(\hat{\l}_r^*,\hat{\S})$ is the same as that of the optimal ridge estimator, $\hat{\bb}_r^*$.  It is worth pointing out that, under the asymptotic setting described in Proposition 11, $n^{-1}X^TX$ is not operator norm-consistent for $\S$; indeed, $\tilde{\bb}_r(\hat{\l}_r^*,n^{-1}X^TX) = \check{\bb}_{js}$ is the adaptive James-Stein estimator and Proposition 10 implies that $\check{\bb}_{js}$ is not asymptotically equivalent to the oracle ridge estimator.  On the other hand, norm-consistent estimators for $\S$ may be available over wide classes of covariance matrices, subject to certain restrictions \citep{bickel2008regularized,elkaroui2008operator,cai2010optimal}.    

\subsection{The James-Stein estimator: Baranchik's estimator}

\citet{baranchik1973inadmissibility} studied the predictive risk of a James-Stein type estimator different from our adaptive James-Stein estimator, $\check{\bb}_{js} = \hat{\bb}(\hat{\l}_{js}^*)$.  Baranchik proved that the estimator 
\begin{equation}\label{bar}
\hat{\bb}_{bar} = \hat{\bb}_{js}(\hat{\l}_{bar}) = (1 + \hat{\l}_{bar})^{-1}\hat{\bb}_{ols} = \left(1 -  c\frac{||\y - X\hat{\bb}_{ols}||^2}{||X\hat{\bb}_{ols}||^2}\right)\hat{\bb}_{ols}
\end{equation}
has smaller predictive risk than the OLS estimator (thus, is minimax over $V \in PD(d+1)$) whenever $d \geq 3$ and $n - d\geq  2$, provided the constant $c$ satisfies $0 < c < 2(d-2)/(n - d + 2)$.  

It is informative to consider the asymptotic predictive risk of Baranchik's estimator.  First notice that (\ref{bar}) implies
\[
\hat{\l}_{bar} = \frac{c ||\y - X\hat{\bb}_{ols}||^2}{||X\hat{\bb}_{ols}||^2 - c||\y - X\hat{\bb}_{ols}||^2}.
\]
The key observation is that if $d, n$ are large and $d/n$ is close to $\rho \in (0,1)$, then
\[
 \hat{\l}_{bar} \approx \l_{bar} \defn \frac{c(1 - \rho)}{\eta^2 + \rho - c(1 -  \rho)},
\]
and, furthermore, $\l_{bar}$ is {\em not} in general equal to the limiting optimal shrinkage parameter,
\[
\l_{js} \defn \lim_{d/n \to \rho} \l_{js}^* = \frac{\rho}{\eta^2(1 - \rho)}.
\]
This suggests that the asymptotic predictive risk of $\hat{\bb}_{bar}$ is suboptimal.  Carrying this heuristic a step further, we take the limit as $n \to \infty$ and $d/n \to \rho \in (0,1)$, and utilize Proposition 4 to obtain an expression for the asymptotic predictive risk of James-Stein type estimators with an arbitrary shrinkage parameter $\l \geq 0$:
\[
R_{js}(\l, \rho,\eta^2) \defn \lim_{d/n \to \rho} R_V\{\hat{\bb}_{js}(\l)\} = \frac{\l^2\eta^2 + \rho/(1-\rho)}{(1 + \l)^2}.
\]
This yields the approximate inequality
\[
R_V(\hat{\bb}_{bar}) \approx  R_{js}(\l_{bar},\rho,\eta^2) 
\geq \min_{\l \geq 0} R_{js}(\l,\rho,\eta^2) =   R_{js}(\l_{js},\rho,\eta^2) \approx R_V(\check{\bb}_{js}),
\]
where $\check{\bb}_{js}$ is the adaptive James-Stein estimator defined in Proposition 10 and the approximation is valid for large $n$ and $d/n$ close to $\rho$. It should be noted that the inequality $R_{js}(\l_{bar}, \rho,\eta^2) \geq R_{js}(\l_{js},\rho,\eta^2)$ is strict unless $\l_{bar} = \l_{js}$.   Moreover, though the equality $\l_{bar} = \l_{js}$ may hold for some specific values of $c$, $\rho$, and $\eta^2$, in order for it to hold in general, the constant $c$ from Baranchik's estimator must vary with $\rho$ and $\eta^2$.  

Some of the ideas from the previous discussion are made more rigorous in the next proposition, whose proof is omitted (part (a) is a straightforward calculation, the proof of (b) is similar to that of Proposition 10, and (c) follows directly from part (b) and Proposition 10).  

\vspace{.1in}
{\em Proposition 12.}
\begin{itemize}
\item[(a)] Suppose $\rho \in (0,1)$.  Then 
\[
 R_{js}(\l_{js},\rho,\eta^2) \leq  R_{js}(\l_{bar},\rho,\eta^2)
\]
with equality if and only if 
\[
c = \frac{\eta^2\rho + \rho^2}{\eta^2(1 - \rho)^2 + \rho(1-\rho)}.
\]
\item[(b)] Suppose that $0 < \th \leq d/n \leq \Th < 1$ for some fixed constants $\th,\Th \in \R$ and that $c$ is a positive constant satisfying $0 < c <  2(d-2)/(n-d+2)$ for all $n$ and $d$.  Further suppose that $n-d > 6$ and let $\rho = d/n$.  Then
\[
R_V(\hat{\bb}_{bar}) = R_V\{\hat{\bb}_{js}(\l_{bar})\} + O\left\{\frac{1}{(\eta^2 + 1)n^{1/2}}\right\}.
\]
\item[(c)]  Under the assumptions of part (b),
\[
R_V(\hat{\bb}_{bar}) - R_V(\check{\bb}_{js})=  R_{js}(\l_{bar},d/n,\eta^2) - R_{js}(\l_{js},d/n,\eta^2)  + O\left\{\frac{1}{(\eta^2 + 1)n^{1/2}}\right\},
\]
where $\check{\bb}_{js}$ is the adaptive James-Stein estimator from Proposition 10. \hfill $\Box$
\end{itemize}
\vspace{.1in}

Proposition 12 and the preceding discussion imply that $\hat{\bb}_{bar}$ is suboptimal in terms of predictive risk, even among the class of James-Stein estimators, $\hat{\bb}_{js}(\l)$.  This naturally leads to the question:  are there other circumstances under which Baranchik's estimator is asymptotically optimal among James-Stein estimators?  The answer is affirmative.  A straightforward calculation (details omitted) implies that $\hat{\bb}_{bar}$ is asymptotically optimal among James-Stein estimator for {\em in-sample} prediction, where the relevant risk function evaluated at an estimator $\hat{\bb}$ is given by
\[
\frac{1}{\s^2n}E_V||X(\hat{\bb} - \bb)||^2.
\]

\subsection{The marginal estimator: Shrinkage}

In Section 5.5 and Proposition 7 we argued that the marginal estimator $\hat{\bb}_m$ has significant drawbacks for out-of-sample prediction.  One natural modification of the marginal estimator that could address some of these drawbacks is the marginal shrinkage estimator,
\[
\hat{\bb}_m(\l) = \frac{1}{(1 + \l)n}\S X^T\y, \ \ \l \geq 0.  
\]
Proposition 13 summarizes some properties of the marginal shrinkage estimator that are analogous to properties of the James-Stein and ridge estimators studied above.

\vspace{.1in}

{\em Proposition 13.} 
\begin{itemize}
\item[(a)] {\em Oracle estimator.} Let $\l_m^* = d/(n\eta^2) + (d + 1)/n$ and let $\hat{\bb}_m^* = \hat{\bb}_m(\l_m^*)$.  Then
\[
R_V(\hat{\bb}_m^*) = \inf_{\l \in [0,\infty]} R_V\{\hat{\bb}_m(\l)\} = \frac{\eta^2\{\eta^2(d+1) + d\}}{\eta^2(n + d + 1) + d}.
\]
\item[(b)] {\em Adaptive estimator.}  Suppose that $0 < \th \leq d/n \leq \Th < 1$ for some fixed $\th,\Th \in \R$ and that $n-d > 4$.  Let $\hat{\l}_m^* = d/(n\hat{\eta^2}) + (d+ 1)/n$ and define the adaptive shrinkage estimator $\check{\bb}_m = \hat{\bb}_m(\hat{\l}_m^*)$.  Then
\[
R_V(\check{\bb}_m) = R_V(\hat{\bb}_m^*) + O(n^{-1/2}).
\]\hfill $\Box$.  
\end{itemize}

The proof of Proposition 13 is omitted.  It may be proved using the same techniques used to prove the analogous results for the ridge and James-Stein estimators. Proposition 13 suggests that the marginal shrinkage estimator still performs relatively poorly when $\eta^2$ is large.  Indeed, 
\begin{equation}\label{margrat}
\lim_{\eta^2 \to \infty} \frac{R_V(\hat{\bb}_m^*)}{R_V(\hat{\bb}_{ols})} > 1 
\end{equation}
In fact, the limit (\ref{margrat}) is infinite if $d < n-1$.  Contrast (\ref{margrat}) with (\ref{ratlim}), which states that the corresponding limit for the oracle ridge and James-Stein estimators equals 1.  Additionally, note that the error term in part (c) of Proposition 13 is independent of $\eta$, while the corresponding error terms in Proposition 10 for the adaptive ridge and James-Stein estimators are proportional to $(1 + \eta^2)^{-1}$. 

Proposition 12 (a) implies that the asymptotic predictive risk of $\hat{\bb}_m^*$ as $d/n \to \rho$ is
\[
R^*_m(\rho,\eta^2) = \frac{\eta^2(\eta^2 + 1)\rho}{\eta^2(1 + \rho) + \rho}.
\]
It is elementary (though somewhat tedious) to check that 
\[
R_r(\rho,\eta^2) \leq R_m^*(\rho,\eta^2)
\]  
for all $\rho$ and $\eta^2$, with equality if and only if $\eta^2 = 0$. In other words, the oracle ridge estimator asymptotically dominates the oracle marginal shrinkage estimator.   On the other hand, neither the oracle marginal shrinkage estimator nor the oracle James-Stein estimator dominates the other.  This is made clear in the plots of asymptotic predictive risk found in Figures 2-3.

\section{Discussion}

\subsection{Adaptive estimators for $d \geq n$}

One limitation of the adaptive estimators considered in Sections 7 and 8 of this paper is the requirement $d < n$.  This requirement is related to the fact that if $d \geq n$, then the usual estimator for $\s^2$, $\hat{\s}^2 = (n-d)^{-1}||\y - X\bb_{ols}||^2$, is undefined.  By replacing $\hat{\s}^2$ in (\ref{esnr}) with an alternative estimator for $\s^2$ that is defined for all $d$, $n$, one could obtain an estimator for $\eta^2$ that is defined for all $d$, $n$.  This would immediately yield adaptive ridge and James-Stein estimators for all $d,n$.  However, estimating $\s^2$ in settings where $d > n$ is challenging.  Recent work by \cite{fan2012variance} and \cite{sun2011scaled} has considered estimating $\s^2$ in settings where $\sup d/n = \infty$, provided $\bb$ is sparse.  We conjecture that even if $\bb$ does not satisfy any sparsity conditions, it may be possible to effectively estimate $\s^2$ when $d \geq n$, provided $\sup d/n < \infty$, and that a adaptive ridge or James-Stein estimator based on such an estimate of $\s^2$ may have the same asymptotic predictive risk as the corresponding oracle estimator.  

\subsection{Conclusions}

Motivated by questions about non-sparse signals in high-dimensional data
analysis, we studied the predictive risk of the OLS, James-Stein, ridge, and marginal regression estimators in high-dimensional linear models.  Our analysis provides new, practical insights into the performance of these popular methods,  while making no assumptions about sparsity.  Both the ridge and James-Stein estimators may substantially outperform the OLS estimator in terms of predictive risk, especially if the signal-to-noise ratio is small.  Additionally, the ridge estimator studied here leverages the population predictor covariance matrix $\S$ to obtain further improvements in out-of-sample prediction when compared to the James-Stein estimator.  These improvements may be precisely quantified by our formulas for asymptotic predictive risk.  On the other hand, we also showed that the marginal regression estimator has substantial drawbacks for out-of-sample prediction.

\section*{Acknowledgments} The author thanks Xihong Lin, Bill Strawderman, Cun-Hui Zhang, and Sihai Zhao for reading various versions of the manuscript and providing helpful comments.  

\section*{Appendices}

\subsection*{Appendix A: Required lemmas}

\vspace{.1in}

{\em Lemma A1. } Let $s_d \geq 0$ denote the smallest eigenvalue of $n^{-1}X^TX$.  Suppose that $k > 0$ is fixed and that $n-d > 2k + 1$. Then
\begin{equation}\label{lemmaA1statement}
E_I(s_d^{-k}) \leq 2\left\{\sqrt{\frac{nd}{8\pi(n-d)}}e^{n+2}\right\}^{2k/(n-d-1)}.
\end{equation}

\vspace{.1in}

{\em Proof.} If $d = 1$, then it is easy to check that the lemma is true.  Thus, suppose that $d \geq 2$.  \citet{muirhead1982aspects} gives the joint density of the ordered eigenvalues, $s_1 > \cdots > s_d > 0$, of $n^{-1} X^TX$:
\[
c_{n,d}\exp\left(-\frac{n}{2}\sum_{j = 1}^d s_j\right)\prod_{j = 1}^d s_j^{(n - d - 1)/2}\prod_{i < j}(s_i - s_j),
\]
where 
\[
c_{n,d} = \frac{\pi^{d^2/2}}{(2/n)^{dn/2}\Gamma_d(d/2)\Gamma_d(n/2)}
\]
and
\[
\Gamma_d(n/2) = \pi^{d(d-1)/4}\prod_{j = 1}^d \Gamma\left\{(n - j + 1)/2\right\}
\]
is the multivariate gamma function.  From this it follows that
\begin{equation}\label{lemmaA1a}
E_I(s_d^{-k}) \leq  E_I(s_d^{-k}; s_d \leq c) + c^{-k} \leq  c^{(n - d - 1)/2 - k}\frac{c_{n,d}}{c_{n,d-1}}E\left\{\det(n^{-1}Z^TZ)^{1/2}\right\} + c^{-k} 
\end{equation}
for any $c > 0$, where $Z$ is an $n \times (d-1)$ dimensional matrix of iid standard normal random variables.   It is easy to check that
\[
\frac{c_{n,d}}{c_{n,d-1}}  =  \frac{\sqrt{\pi}(n/2)^{n/2}}{\Gamma\{(n - d + 1)/2\}\Gamma(d/2)}.
\]
Additionally, it is well known (Exercise 3.11 in \citet{muirhead1982aspects}, for instance) that
\[
E\left\{\det(n^{-1}Z^TZ)^{1/2}\right\}  =  (2/n)^{(d-1)/2}\frac{\Gamma_{d-1}\{(n + 1)/2\}}{\Gamma_{d-1}(n/2)}  =  (2/n)^{(d-1)/2}\frac{\Gamma\{(n + 1)/2\}}{\Gamma\{(n - d + 2)/2\}}
\]
Thus, using Stirling's approximation, we obtain
\begin{eqnarray*}
\frac{c_{n,d}}{c_{n,d-1}} E\left\{\det(n^{-1}Z^TZ)^{1/2}\right\} & = & \frac{\sqrt{\pi}(n/2)^{(n-d+1)/2}\Gamma\{(n + 1)/2\}}{\Gamma\{(n - d + 1)/2\} \Gamma\{(n - d + 2)/2\} \Gamma(d/2)}\\
& = & \frac{(2n)^{(n-d+1)/2}\Gamma\{(n + 1)/2\}}{2\Gamma(n-d+1)\Gamma(d/2)} \\
& \leq & \sqrt{\frac{nd}{8\pi(n-d)}}  \left(\frac{n}{n-d}\right)^{n-d} \left(\frac{n}{d}\right)^{d/2}e^{(n-d)/2+2} \\
& \leq & \sqrt{\frac{nd}{8\pi(n-d)}}e^{n+2}
\end{eqnarray*}
and, by (\ref{lemmaA1a}),
\[
E_I(s_d^{-k})  \leq   c^{(n-d-1)/2 - k} \sqrt{\frac{nd}{8\pi(n-d)}}e^{n+2} + c^{-k}
\]
Taking $c =  \left\{\sqrt{nd/\{8\pi(n-d)\}}e^{n+2}\right\}^{-2/(n-d-1)}$ gives (\ref{lemmaA1statement}).  
\hfill $\Box$

\vspace{.1in}

{\em Lemma A2.} Let $s_1 \geq s_d \geq 0$ denote the largest and smallest eigenvalues of $n^{-1}X^TX$, respectively.  Suppose that $k > 0$ is fixed and that $0 < d/n \leq \Th < 1$ for some fixed constant $\Th \in \R$.  
\begin{itemize}
\item[(a)] $E_I(s_1^k) = O(1)$.
\item[(b)] If $n-d > 2k + 1$, then $E_I(s_d^{-k}) = O(1)$.  
\end{itemize}

\vspace{.1in}

{\em Proof.} Part (a) is well known and may be easily derived from large deviations results for $s_1$ (see, for example, Theorem II.13 of \citep{davidson2001local}). Part (b) follows directly from Lemma A2.  
\hfill $\Box$

\vspace{.1in}

{\em Lemma A3.} Let $k$ be a fixed positive integer and suppose that $0  <  d/n \leq \Th < 1$ for some fixed constant $ \Th \in \R$.  Then
\[
E_V\left(\frac{1}{\hat{\eta}^2 + d/n}\right)^k= O\left\{\left(\frac{1}{\eta^2 + d/n}\right)^{k}\right\}.
\]

\vspace{.1in}

{\em Proof.} 
Notice that
\begin{eqnarray*}
E_V\left(\frac{1}{\hat{\eta}^2 + 2d/n}\right)^k & \leq & E_V\left\{\frac{1}{||X\hat{\bb}_{ols}||^2/(n\hat{\s}^2) + d/n}\right\}^k \\
& \leq& E_V\left\{\frac{\hat{\s}^2/\s^2 + 1}{||X\hat{\bb}_{ols}||^2/(n\s^2) + d/n}\right\}^k\\
& = &  E_V\left(\frac{\hat{\s}^2}{\s^2} + 1 \right)^k E_V\left\{\frac{1}{||X\hat{\bb}_{ols}||^2/(n\s^2) + d/n}\right\}^k
\end{eqnarray*}
and that $E_V(\hat{\s}^2/\s^2 + 1)^k = O(1)$.  Since it is clear that $E_V\left(\hat{\eta}^2 + d/n\right)^{-k}=O(n/d)$,  it suffices to show that 
\[
E_V\left\{\frac{1}{||X\hat{\bb}_{ols}||^2/(n\s^2) + d/n}\right\}^k= O(\eta^{-2k}).
\]
Conditional on $X$, $||X\hat{\bb}_{ols}||^2/\s^2$ follows a noncentral $\chi^2$ distribution with $d$ degrees of freedom and noncentrality parameter $||X\bb||^2/\s^2$.  Thus $||X\hat{\bb}_{ols}||^2/\s^2$ has the same distribution as a central $\chi^2$ random variable with $2N + d$ degrees of freedom, where $N|X \sim \mbox{Poisson}(\zeta)$ and $\zeta = ||X\bb||^2/(2\s^2)$.  Now let $W \sim \chi^2_m$ be an independent $\chi^2$ random variable with $m = (2k-d+2)\vee 1$ degrees of freedom.  Then, since the $k$-th inverse moment of a $\chi^2$ random variable with $l$ degrees of freedom is $2^{-k} \Gamma(l/2-k)/\Gamma(l/2)$, provided $l > 2k$, it follows from Jensen's inequality that
\begin{eqnarray}
E_V\left\{\frac{1}{||X\hat{\bb}_{ols}||^2/(n\s^2) + d/n}\right\}^k & \leq & E_V\left\{\frac{1}{||X\hat{\bb}_{ols}||^2/(n\s^2) +  d/(mn)W}\right\}^k \nonumber \\
& \leq &  \left\{\frac{n(m+d)}{d}\right\}^k  E_V\left(\frac{1}{||X\hat{\bb}_{ols}||^2/\s^2 +W}\right)^k  \nonumber \\
& = &   \left\{\frac{n(m+d)}{2d}\right\}^k E\left[\frac{\Gamma\{N + (d+m)/2 - k\}}{\Gamma\{N + (d+m)/2\}}\right] \nonumber \\
& = &  \left\{\frac{n(m+d)}{2d}\right\}^k E\left\{\prod_{i = 1}^k\frac{1}{N + (d+m)/2 - i}\right\} \nonumber \\
& \leq &   \left\{\frac{n(m+d)}{2d}\right\}^k\sum_{j = 0}^{\infty} E_V\left\{\frac{\zeta^j}{(j + k)!}e^{-\zeta}\right\} \nonumber \\
& \leq &   \left\{\frac{n(m+d)}{2d}\right\}^kE_V(\zeta^{-k}) \nonumber \\
& = & O(\eta^{-k}),
\end{eqnarray}
where we have used the fact that $\zeta/\eta^2 \sim \chi^2_n$ has a $\chi^2$ distribution with $n$ degrees of freedom.  The lemma follows.   \hfill $\Box$

\vspace{.1in}

{\em Lemma A4. } Suppose that $0 < d/n \leq \Theta < 1$ for some fixed constant $\Th \in \R$ and let $k > 0$ be fixed.  If $n > 2k$, then 
\[
 P_V(\hat{\eta}^2 = 0) = O\left(\frac{d^{k/2}}{\eta^{2k}n^{k}}\right). 
\]

\vspace{.1in} 

{\em Proof.}  Let $U = ||X\hat{\bb}_{ols}||^2/\s^2$ and let $W = ||y - X\hat{\bb}_{ols}||^2/\s^2 = (n-p)\hat{\s}^2/\s^2$.   Then $W \sim \chi^2_{n-d}$ has a $\chi^2$ distribution with $n-d$ degrees of freedom and, conditional on $X$, $U\sim \chi^2_{||X\bb||^2/\s^2,d}$ has a noncentral $\chi^2$ distribution with noncentrality parameter $||X\bb||^2/\s^2$ and $d$ degrees of freedom.  Furthermore, $U$ and $W$ are independent and
\[
\hat{\eta}^2 = \frac{d}{n}\left\{\frac{U/d}{W/(n-d)} - 1\right\} \vee 0.
\]
Thus,
\begin{eqnarray*}
P_V(\hat{\eta}^2 = 0) & = & P_V\left(\frac{1}{d}U \leq \frac{1}{n-d}W\right) \\
& \leq & E_V\exp\left(\frac{r}{n-d}W\right)E_V\left(-\frac{r}{d}U\right) \\
& = & \left(\frac{1}{1 - \frac{2r}{n-d}}\right)^{(n-d)/2} \left(\frac{1}{1 + \frac{2r}{d}}\right)^{d/2} E_V\exp\left(-\frac{r}{d + 2r}||X\bb||^2/\s^2\right) \\
& = & \left(\frac{n-d}{n-d - 2r}\right)^{(n-d)/2}\left(\frac{d}{d + 2r}\right)^{d/2} \left\{\frac{d + 2r}{d + 2(\eta^2 + 1)r}\right\}^{n/2} \\
& \leq & \exp\left\{\frac{2r^2n}{(n-d-2r)(d + 2r)}\right\}\left\{\frac{d + 2r}{d + 2(\eta^2 + 1)r}\right\}^{n/2},
\end{eqnarray*}
provided $r < (n-d)/2$.  Now, basic calculus implies that 
\[
\sup_{\eta^2 \geq 0} \eta^{2k}\left\{\frac{d + 2r}{d + 2(\eta^2 + 1)r}\right\}^{n/2} \leq e^{-k}\left\{\frac{(d + 2r)k}{r(n - 2k)}\right\}^k.  
\]
The lemma follows by taking $r = \a \sqrt{d}$ for $\a > 0$ sufficiently small.  
\hfill $\Box$

\vspace{.1in}

{\em Lemma A5. }
Suppose that $0 < d/n \leq\Theta < 1$ for some fixed constant $\Theta \in \R$ and that $k > 0$ is fixed.  If $n - d > 2k$, then 
\[
E_V|\hat{\eta}^2 - \eta^2|^k = O\left\{\left(\frac{d/n + \eta^2}{n}\right)^{k/2} + \frac{\eta^{2k}}{n^{k/2}}\right\}.
\]

\vspace{.1in}

{\em Proof.} Using Lemma A4, we have
\begin{eqnarray}
E_V|\hat{\eta}^2 - \eta^2|^k & \leq&  E_V\left|\frac{||X\hat{\bb}_{ols}||^2}{n\hat{\s}^2} - \left(\frac{d}{n} + \eta^2\right)\right|^k + \eta^{2k}P_V(\hat{\eta}^2 = 0) \nonumber \\
& = &  E_V\left|\frac{||X\hat{\bb}_{ols}||^2}{n\hat{\s}^2}  - \left(\frac{d}{n} + \eta^2\right)\right|^k  + O\left(\frac{d^{k/2}}{n^k}\right). \label{lemmaA5a}
\end{eqnarray}
Since $(n-d)\hat{\s}^2 = ||\y - X\hat{\bb}_{ols}||^2$ and $||X\hat{\bb}_{ols}||^2$ are independent,
\begin{eqnarray*}
E_V\left|\frac{||X\hat{\bb}_{ols}||^2}{n\hat{\s}^2} - \left(\frac{d}{n} + \eta^2\right)\right|^k & \leq & 2^kE_V\left|\frac{||X\hat{\bb}_{ols}||^2}{n\hat{\s}^2} - \left(\eta^2 + \frac{d}{n}\right)\frac{\s^2}{\hat{\s}^2}\right|^k \\
&& \quad + 2^k \left(\eta^2 + \frac{d}{n}\right)^k E_V\left|\frac{\s^2}{\hat{\s}^2} - 1\right|^k \\
& = & 2^kE_V\left(\frac{\s^2}{\hat{\s}^2}\right)^k E_V\left|\frac{||X\hat{\bb}_{ols}||^2}{n\s^2} - \left(\eta^2 + \frac{d}{n}\right)\right|^k \\
&& \quad + 2^k \left(\eta^2 + \frac{d}{n}\right)^k E_V\left|\frac{\s^2}{\hat{\s}^2} - 1\right|^k. 
\end{eqnarray*}
As in the proof of Lemma A3, let $N \sim \mbox{Poisson}\{||X\bb||^2/(2\s^2)\}$.  Since $||X\hat{\bb}_{ols}||^2/\s^2 \sim \chi^2_{2N + d}$, we have  
\begin{eqnarray*}
E_V\left|\frac{||X\hat{\bb}_{ols}||^2}{n\s^2} - \left(\eta^2 + \frac{d}{n}\right)\right|^k & \leq & 2^kE_V\left|\frac{||X\hat{\bb}_{ols}||^2}{n\s^2} - \frac{2N+d}{n}\right|^k +2^k E_V\left|\frac{2N}{n} - \eta^2\right|^k \\
& = & O\left\{n^{-k/2}\left(\eta^2 +\frac{d}{n}\right)^{k/2}\right\}
\end{eqnarray*}  
Additionally, one can check that 
\[
 E_V\left|\frac{\s^2}{\hat{\s}^2} - 1\right|^k = n^{-k/2}.
\]
It follows that
\[
E_V\left|\frac{||X\hat{\bb}_{ols}||^2}{n\hat{\s}^2} - \left(\frac{d}{n} + \eta^2\right)\right|^k = O\left\{\frac{(d/n)^{k/2} + \eta^k + \eta^{2k}}{n^{k/2}}\right\}.
\]
The lemma follows by combining this with (\ref{lemmaA5a}).  
 \hfill $\Box$

\subsection*{Appendix B: Proofs of propositions contained in the main text}\hspace{.1in} \newline

{\em Proof of Proposition 1.}  Suppose that $\hat{\bb}$ is linearly equivariant.  To prove part (a), observe that
\begin{eqnarray*}
\s^2R_V(\hat{\bb}) & = & E_V\left\{\hat{\bb}(y,X,\S) - \bb\right\}^T\S\left\{\hat{\bb}(y,X,\S) - \bb\right\} \\ &  = & E_V\left|\left|\hat{\bb}(y,X\S^{-1/2},I) - \S^{1/2}\bb\right|\right|^2 \\
& = & E_{V_0}\left|\left|\hat{\bb}(y,X,I) - \bb(V_0)\right|\right|^2 \\
& = & \s^2R_{V_0}(\hat{\bb}),
\end{eqnarray*}
where $\bb(V_0) = \S^{1/2}\bb$, and we have used linear equivariance of $\hat{\bb}$, along with the fact that if $\x_i \sim N(0,\S)$, then $\S^{-1/2}\x_i \sim N(0,I)$.  

Now suppose that $\hat{\bb}$ is LS, let $\u \in \R^d$ be a unit vector, and let $U$ be a $d \times d$ orthogonal matrix such that $U\S^{1/2}\bb/\s = \eta \u \defn \bb(V_{\u})$.  Then
\begin{eqnarray*}
R_V(\hat{\bb}) & = & \s^{-2}E_{V_0}\left|\left|\hat{\bb} - \S^{1/2}\bb\right|\right|^2 \\
& = & E_{V_0}\left|\left|\s^{-1}U\hat{\bb}(y,X,I) - \eta\u\right|\right|^2 \\
& = &  E_{V_0}\left|\left|\hat{\bb}(XU^T \eta\u + \ee/\s, XU^T,I) - \eta\u\right|\right|^2 \\
& = & E_{V_{\u}}\left|\left|\hat{\bb}(y,X,I) - \bb(V_{\u})\right|\right|^2 \\
& = & R_{V_{\u}}(\hat{\bb}).
\end{eqnarray*}
\hfill $\Box$

\vspace{.1in}

{\em Proof of Proposition 2.} Proposition 1 and a simple bias-variance decomposition lead to 
\[
R_V(\hat{\bb}_{ols}) = R_{V_{\u}}(\hat{\bb}_{ols}) = \eta^2 E_I\left|\left|\left\{(X^TX)^-X^TX-I\right\}\u\right|\right|^2 +  E_I\tr\left\{(X^TX)^-\right\} 
\]
If $d \leq n$, then the first term on the right-hand side above is equal to 0; if $d > n$, then, by symmetry, it is equal to $\eta^2(d-n)/d$.  Using properties of the inverse Wishart distribution (see Chapter 3 of \citep{muirhead1982aspects}, for instance) the second term on the right-hand side is equal to $d/(n-d-1)$, if $d < n-1$; it is equal to $n/(d-n-1)$ if $d > n+1$; and it is infinite otherwise.  The proposition follows. \hfill $\Box$

\vspace{.1in}

{\em Proof of Proposition 3.} Fix a unit vector $\u = (u_1,...,u_d)^T \in \R^d$.  Then
\begin{equation}\label{prop2a}
R_V(\hat{\bb}_m) = E_{V_{\u}}||n^{-1}X^T\y - \eta \u||^2 = \sum_{j = 1}^d E_{V_{\u}}\left(\frac{1}{n}\X_j^T\y - \eta u_j\right)^2,
\end{equation}
where $\X_j$ is the $j$-th column of $X$.  Considering each term in the summation above separately, we have
\begin{eqnarray*}
E_{V_{\u}}\left(\frac{1}{n}\X_j^T\y - \eta u_j\right)^2 & = & E_{V_{\u}}\left\{\left(\frac{1}{n}\X_j^T\X_j - 1\right)\eta u_j  + \frac{1}{n}\sum_{k \neq j} \X_j^T\X_k\eta u_k + \frac{1}{n}\X_j^T\ee\right\}^2 \\
& = & \eta^2 u_j^2E_I\left(\frac{1}{n}\X_j^T\X_j - 1\right)^2   \\
&& \qquad + \frac{\eta^2}{n^2} E_I\left(\sum_{k \neq j} \X_j^T\X_k u_k\right)^2 + \frac{1}{n^2}E_{V_{\u}}(\X_j^T\ee)^2 \\
& = & \frac{2\eta^2u_j^2}{n} +\frac{\eta^2}{n}\sum_{k \neq j} u_k^2 + \frac{1}{n} \\
& = & \frac{1}{n}(\eta^2 u_j^2 + \eta^2 + 1).
\end{eqnarray*}
The proposition follows by summing over $j$ above and using the identity (\ref{prop2a}). \hfill $\Box$

\vspace{.1in}

{\em Proof of Proposition 4.}  It is easy to check that
\[
R_{V_{\u}}\{\hat{\bb}_{js}(\l)\}  =  \left(\frac{1}{1 + \l}\right)^2J_1 + \left(\frac{\l}{1 + \l}\right)^2J_2 +J_3.
\]
where
\begin{eqnarray*}
J_1 & = & E_I\tr\left\{(X^TX)^-\right\} \ \ = \ \ \left\{\begin{array}{cl} \frac{d}{n-d-1} & \mbox{if } d < n -1 \\
\frac{n}{d-n-1} & \mbox{if } d > n + 1 \\
\infty & \mbox{if } p \in \{n-1,n,n+1\} \end{array}\right.  \\
J_2 & = & \eta^2E_I\left|\left|(X^TX)^-X^TX\u\right|\right|^2 \ \ = \ \ \left\{\begin{array}{cl} \eta^2 & \mbox{if } d \leq n \\  \eta^2\frac{n}{d} & \mbox{if } d > n \end{array} \right.   \\
J_3 & = & \eta^2 E_I\left|\left|\left\{I - (X^TX)^-X^TX\right\}\bb\right|\right|^2 \ \ = \ \ \left\{\begin{array}{cl} 0 & \mbox{if } d \leq n \\
\frac{d-n}{d}\eta^2 & \mbox{if } d > n.  \end{array}\right.  
\end{eqnarray*}
Hence, (\ref{prop3a}).  The rest of the proposition follows by basic calculus.  \hfill $\Box$
\vspace{.1in}

{\em Proof of Proposition 5.} For $j = 1,...,d$, let $\mathbf{e}_j \in \R^d$ denote the $j$-th standard basis vector.  Fix $\l \in [0,\infty]$.  Then
\begin{eqnarray*}
R_V\{\hat{\bb}_r(\l)\} & = & E_{V_{\mathbf{e}_j}}\left|\left|\hat{\bb}_r(\l) - \bb\right|\right|^2 \\
& = & E_{V_{\mathbf{e}_j}}\left|\left|\eta n\l(X^TX + n\l I)^{-1}\mathbf{e}_j\right|\right|^2  +  E_{V_{\mathbf{e}_j}}\left|\left|(X^TX + n\l I)^{-1}X^T\ee\right|\right|^2 \\
& = &  E_I\left|\left|\eta n\l(X^TX + n\l I)^{-1}\mathbf{e}_j\right|\right|^2 + E_I\tr\left\{(X^TX + n \l I)^{-2}X^TX\right\}.
\end{eqnarray*}
Summing over $j = 1,...,d$ above gives
\[
R_V\{\hat{\bb}_r(\l)\} = \frac{\eta^2n^2}{d}\l^2 E_I\tr\left\{(X^TX + n\l I)^{-2}\right\} + E_I\tr\left\{(X^TX + n \l I)^{-2}X^TX\right\}
\]
and (\ref{prop6a}) follows.  

To prove (\ref{prop6b}), let $s_1 \geq \cdots s_d \geq 0$ denote the eigenvalues of $n^{-1}X^TX$.  Then
\[
R_V\{\hat{\bb}_r(\l)\}  =  E_I\left\{\sum_{j = 1}^d (ns_j + n\l)^{-2}\left(ns_j + \frac{\eta^2n^2}{d}\l\right)\right\}.
\]
It is easy to check that each of the $d$ summands on the right-hand side above is minimized by taking $\l = \l_r^*$ and that $R_V(\hat{\bb}_r^*) = E_I\tr(X^TX + n\l_r^*I)^{-1}$.  \hfill $\Box$  

\vspace{.1in}

{\em Proof of Proposition 6.} The inequality $R_V(\hat{\bb}_{js}^*) < R_V(\hat{\bb}_{ols})$ was discussed above and  follows from Propositions 2 and 4.  The other inequality follows from Jensen's inequality.  As in the proof of Proposition 5, let $s_1 \geq  \cdots \geq s_d \geq 0$ denote the eigenvalues of $n^{-1}X^TX$.  First suppose that $d  < n-1$.  Then
\begin{eqnarray*}
R_V(\hat{\bb}_r^*) & = & E_I\tr(X^TX + n\l_r^*I)^{-1} \\
& = & E_I\left\{\frac{1}{n}\sum_{j = 1}^d \frac{1}{s_j + \l^*_r} \right\} \\
& = & \eta^2E_I\left\{\frac{1}{d} \sum_{j = 1}^d \frac{s_j^{-1}}{n\eta^2/d + s_j^{-1}} \right\} \\
& \leq & \eta^2 \frac{E_I\tr(X^TX)^{-1}}{\eta^2 + E_I\tr(X^TX)^{-1}} \\
& = & \frac{\eta^2d}{\eta^2(n - d - 1) + d} \\
& = & R_V(\hat{\bb}_{js}^*),
\end{eqnarray*}  
where the inequality is strict unless $\eta^2 = 0$.  If $d > n + 1$, then $s_{n+1} = s_{n+2} = \cdots = s_d = 0$ and a similar calculation implies that $R_V(\hat{\bb}_r^*) \leq R_V(\hat{\bb}_{js}^*)$.
Finally, if $d \in \{n-1,n,n+1\}$, then it is clear that $R_V(\hat{\bb}_r^*) \leq\eta^2 =  R_V(\hat{\bb}_{js}^*)$.\hfill $\Box$
\vspace{.1in}

\textit{Proof of Proposition 8. }
Let $\F_{n,d}$ be the empirical cumulative distribution function of the eigenvalues of $n^{-1}X^TX$.   Using integration by parts, for $c \geq 0$, 
\begin{eqnarray}
\frac{n}{d}\tr(X^TX + n\l_r^* I)^{-1} & = &  \int_0^{\infty} \frac{1}{s + \l_r^*} \ d\F_{n,d}(s) \nonumber \\
& = & \int_0^c \frac{1}{s + \l_r^*} \ d\F_{n,d}(s) + \frac{1}{c + \l_r^*}\{1 - \F_{n,d}(c)\} \nonumber \\
&& \qquad - \int_c^{\infty} \frac{1}{(s + \l_r^*)^2}\{1 - \F_{n,d}(s)\}  \ ds. \label{prop6proofa}
\end{eqnarray}
Similarly,
\begin{eqnarray}
m_{d/n}(-\l_r^*) & = & \int_0^c \frac{1}{s + \l_r^*} \ dF_{d/n}(s) + \frac{1}{c + \l_r^*}\{1 - F_{d/n}(s)\} \nonumber \\
 && \qquad - \int_c^{\infty} \frac{1}{(s + \l_r^*)^2}\{1 - F_{d/n}(s)\} \ ds. \label{prop6proofb}
\end{eqnarray}

Now let $\Delta = | R_V(\hat{\bb}_r^*) - (d/n)m_{d/n}(-\l_r^*)|$. Then Theorem 3.1 of \citet{bai1993convergence} (see equation (\ref{bai}) in Section 6.1 above) and the inequalities
\[
\Delta  \leq  \frac{d}{n}\int_0^{\infty} \frac{1}{(s + \l_r^*)^2} \left|E_I\F_{n,d}(s) - F_{d/n}(s)\right| \ ds \leq  \eta^2 \sup_{s \geq 0} \left|E_I\F_{n,d}(s) - F_{d/n}(s)\right| ,
\]
which follow from (\ref{prop6proofa}) and (\ref{prop6proofb}), imply
\[
\Delta = \left\{\begin{array}{cl} O(\eta^2 n^{-1/4}) & \mbox{if } 0 < \th < \Th < 1 \mbox{ or } 1 < \th < \Th < \infty, \\
O(\eta^2 n^{-5/48}) & \mbox{if } 0 < \th < 1 < \Th < \infty. \end{array}\right.
\]
Part (b) of the proposition follows immediately.  

To prove part (a) of the proposition, we show that, in fact, $\Delta = O(n^{-1/4})$ if $0 < \th < \Th < 1$ or $1 < \th < \Th < \infty$.  First suppose that $0 < \th < \Th < 1$.  Then, for $c < (1 - \sqrt{d/n})^2$,
\[
m_{d/n}(-\l_r^*) = \frac{1}{c + \l_r^*} - \int_c^{\infty} \frac{1}{(s + \l_r^*)^2}\{1 - F_{d/n}(s)\}  \ ds
\]
and
\begin{eqnarray*}
\frac{n}{d}\Delta & \leq &  E _I\left\{\int_0^c \frac{1}{s + \l_r^*} \ d\F_{n,d}(s)\right\} \\
&& \qquad + \frac{1}{c + \l_r^*} E_I\F_{n,d}(c) + \left| \int_c^{\infty} \frac{1}{(s + \l_r^*)^2}\left\{E_I\F_{n,d}(s) - F_{d/n}(s)\right\} \ ds \right| \\
& \leq & E_I\left\{\int_0^c s^{-1} \ d\F_{n,d}(s)\right\} + \frac{1}{c + \l_r^*} E_I\F_{n,d}(c) + \frac{1}{c + \l_r^*}\sup_{s \geq c} \left|E_I\F_{n,d}(s) - F_{d/n}(s)\right| \\
& \leq & E_I(s_d^{-1}; \ s_d < c) + \frac{1}{c + \l_r^*}P_I(s_d < c) + \frac{1}{c + \l_r^*}\sup_{s \geq c} \left|E_I\F_{n,d}(s) - F_{d/n}(s)\right| \\
& \leq & \left\{E_I(s_d^{-2})\right\}^{1/2}P_I(s_d < c)^{1/2} + c^{-1}P_I(s_d < c)  + c^{-1}\sup_{s \geq c} \left|E_I\F_{n,d}(s) - F_{d/n}(s)\right|,
\end{eqnarray*}
where $ s_d \geq 0$ is the smallest eigenvalue of $n^{-1}X^TX$.  We bound the first two terms and the last term on right-hand side above separately.  Bounding the first two terms relies on a result of \citet{davidson2001local}.  Their Theorem II.13, which is a consequence of concentration of measure,  implies that
\begin{equation} \label{lemma2a}
P_I(s_d \leq c) \leq \exp\left\{-\frac{n(1 - \sqrt{d/n})^2}{2}\left\{1 - \frac{c^{1/2}}{1 - \sqrt{d/n}}\right\}^2\right\},
\end{equation}
provided $c \leq 1 - \sqrt{d/n}$.  Additionally, Lemma A2 in Appendix A implies that $E_I(s_d^{-2}) = O(1)$ if $n - d > 5$.  Taking $c = (1 - \sqrt{d/n})^2/2$, it follows that
\[
\left\{E_I(s_d^{-2})\right\}^{1/2}P(s_d < c)^{1/2} + c^{-1}P(s_d < c) = O(n^{-1/4})
\] 
(in fact, we can conclude that the quantities on the left above decay exponentially, but this is not required for the current proposition). It now follows from Theorem 3.1 of \citet{bai1993convergence} that $\Delta = O(n^{-1/4})$.  For the case where $1 < \th < \Theta$,  we note that the same argument as above applies to $XX^T$, which has the same nonzero eigenvalues as $X^TX$.  Part (a) of the proposition follows.    \hfill $\Box$

\vspace{.1in}

{\em Proof of Proposition 9.} The first statement is easily verified.  To prove the second statement, Proposition 6 implies that it suffices to show
\[
\liminf_{d/n \to 0 \atop d/(n\eta^2)\to c} \frac{R_V(\hat{\bb}_{r}^*)}{R_V(\hat{\bb}_{js}^*)}  \geq 1,
\]
for $c \in [0,\infty]$.  But this follows from Jensen's inequality, which implies 
\[
\frac{d/n}{1 + d/(n\eta^2)} \leq E\tr\{X^TX + (d/\eta^2)I\}^{-1} = R_V(\hat{\bb}_r^*).  
\]
\hfill$\Box$

\vspace{.1in}

\textit{Proof of Proposition 10. }
Assume that $0 < \th \leq d/n \leq \Th < 1$ for some fixed constants $\th,\Th \in \R$ and that $n - d> 9$.  We prove (a).  The proof of (b) is entirely similar.  Since $\hat{\bb}_r^*$ and $\check{\bb}_r$ are LS,
\begin{eqnarray*}
R_V(\hat{\bb}^*_r) & = & R_{V_{\u}}(\hat{\bb}_r^*) \\
& = & E_{V_{\u}}\left|\left|(X^TX + n\l_r^*I)^{-1}X^T\y - \eta \u\right|\right|^2 \\
& = & \eta^2 \left(\frac{d}{n}\right)^2E_{V_{\u}}\left|\left|\left(\eta^2\frac{1}{n}X^TX + \frac{d}{n}I\right)^{-1}\u\right|\right|^2 \\
&& \qquad  - 2\eta^3\frac{d}{n^2}E_{V_{\u}}\ee^TX\left(\eta^2\frac{1}{n}X^TX + \frac{d}{n}I\right)^{-2}\u \\
&& \qquad+\frac{1}{n^2}E_{V_{\u}}\left|\left|\eta^2\left(\eta^2\frac{1}{n}(X^TX + \frac{d}{n}I\right)^{-1}X^T\ee\right|\right|^2  \\
& = & \eta^2 \left(\frac{d}{n}\right)^2E_{V_{\u}}\left|\left|\left(\eta^2\frac{1}{n}X^TX + \frac{d}{n}I\right)^{-1}\u\right|\right|^2 \\
&& \qquad + \frac{1}{n^2}E_{V_{\u}}\left|\left|\eta^2\left(\eta^2\frac{1}{n}X^TX + \frac{d}{n}I\right)^{-1}X^T\ee\right|\right|^2
\end{eqnarray*}
and
\begin{eqnarray*}
R_V(\check{\bb}_r) & = & \eta^2 \left(\frac{d}{n}\right)^2E_{V_{\u}}\left|\left|\left(\hat{\eta}^2\frac{1}{n}X^TX + \frac{d}{n}I\right)^{-1}\u\right|\right|^2 \\
&& \qquad - 2\eta\frac{d}{n^2}E_{V_{\u}}\hat{\eta}^2\ee^TX\left(\hat{\eta}^2\frac{1}{n}X^TX + \frac{d}{n}I\right)^{-2}\u \\
&& \qquad + \frac{1}{n^2}E_{V_{\u}}\left|\left|\hat{\eta}^2\left(\hat{\eta}^2\frac{1}{n}X^TX + \frac{d}{n}I\right)^{-1}X^T\ee\right|\right|^2,
\end{eqnarray*}
where $\u \in \R^d$ is a fixed unit vector.  Thus,
\begin{equation}\label{prop11proofa}
|R_V(\check{\bb}_r) - R_V(\hat{\bb}_r^*)| \leq \left|E_{V_{\u}}H_1\right| + \left|E_{V_{\u}}H_2\right| + 2\left|E_{V_{\u}}H_3\right|, 
\end{equation}
where
\begin{eqnarray*}
H_1 & = & \eta^2\frac{d}{n}\left\{\left|\left|\left(\hat{\eta}^2\frac{1}{n}X^TX + \frac{d}{n}I\right)^{-1}\u\right|\right|^2 - \left|\left|\left(\eta^2\frac{1}{n}X^TX + \frac{d}{n}I\right)^{-1}\u\right|\right|^2\right\} \\
H_2 & = & \frac{1}{n^2}\left\{\left|\left|\hat{\eta}^2\left(\hat{\eta}^2\frac{1}{n}X^TX + \frac{d}{n}I\right)^{-1}X^T\ee\right|\right|^2 - \left|\left|\eta^2\left(\eta^2\frac{1}{n}X^TX + \frac{d}{n}I\right)^{-1}X^T\ee\right|\right|^2\right\} \\
H_3 & = & \eta\hat{\eta}^2\frac{d}{n^2}\ee^TX\left(\hat{\eta}^2\frac{1}{n}X^TX + \frac{d}{n}I\right)^{-2}\u.
\end{eqnarray*}
We consider the terms $|E_{V_{\u}}H_1|$, $|E_{V_{\u}}H_2|$, and $|E_{V_{\u}}H_3|$ separately.

Let $s_1 \geq \cdots \geq s_d \geq 0$ denote the ordered eigenvalues of $n^{-1}X^TX$ and let $U$ be a $d \times d$ orthogonal matrix such that $S = n^{-1}U^TX^TXU$ is diagonal. Additionally, let $\tilde{\u} = (\tilde{u}_1,...,\tilde{u}_d)^T = U^T\u$ and let $\tilde{\dd} = (\tilde{\d}_1,...,\tilde{\d}_d)^T= U^T(X^TX)^{-1/2}X^T\e$.  Then
\begin{eqnarray*}
|H_1 |
&= & \eta^2\frac{d}{n}\left|\sum_{j = 1}^d \left\{\frac{\tilde{u}_j^2}{(\hat{\eta}^2s_j + d/n)^2} -\frac{\tilde{u}_j^2}{(\eta^2s_j + d/n)^2} \right\}\right| \\
& = & \eta^2\frac{d}{n}\left|\sum_{j = 1}^d \frac{\tilde{u}_j^2s_j(\eta^2 - \hat{\eta}^2)}{(\hat{\eta}^2s_j + d/n)(\eta^2s_j + d/n)}\left(\frac{1}{\hat{\eta}^2s_j + d/n} + \frac{1}{\eta^2s_j + d/n}\right)\right| \\
& \leq & 8\eta^2\frac{d}{n}\sum_{j = 1}^d \frac{\tilde{u}_j^2|\eta^2 - \hat{\eta}^2|}{(\hat{\eta}^2 + d/n)(\eta^2 + d/n)}\left(\frac{1}{\hat{\eta}^2 + d/n} + \frac{1}{\eta^2 + d/n}\right)\left(\frac{1}{s_j^{2}} + s_j\right) \\
& \leq & \frac{8(d/n)|\eta^2 - \hat{\eta}^2|}{\hat{\eta}^2 + d/n}\left(\frac{1}{\hat{\eta}^2 + d/n} + \frac{1}{\eta^2 + d/n}\right)\left(\frac{1}{s_d^2} + s_1\right)
\end{eqnarray*}
and 
\begin{eqnarray*}
|H_2| & = & \frac{1}{n} \left|\sum_{j = 1}^d \left\{ \frac{\hat{\eta}^4s_j\tilde{\d}_j^2}{(\hat{\eta}^2s_j + d/n)^2} - \frac{\eta^4s_j\tilde{\d}_j^2}{(\eta^2s_j + d/n)^2}\right\}\right| \\
& = & \frac{1}{n} \left|\sum_{j = 1}^d  \frac{(d/n)\tilde{\d}_j^2s_j(\hat{\eta}^2 - \eta^2)}{(\hat{\eta}^2s_j + d/n)(\eta^2s_j + d/n)}\left(\frac{\hat{\eta}^2}{\hat{\eta}^2s_j + d/n} + \frac{\eta^2}{\eta^2s_j + d/n}\right) \right| \\
& \leq & \frac{4}{n}\sum_{j = 1}^d \frac{\tilde{\d}_j^2|\hat{\eta}^2 - \eta^2|}{(\hat{\eta}^2 + d/n)(\eta^2 + d/n)}\left(\frac{1}{s_j} + s_j\right) \\
& \leq & \frac{4}{n}||\tilde{\dd}||^2\frac{|\eta^2 - \hat{\eta}^2|}{(\hat{\eta}^2 + d/n)(\eta^2 + d/n)}\left(\frac{1}{s_d} + s_1\right).
\end{eqnarray*}
Thus, by Lemmas A2, A3, A5, and the Cauchy-Schwarz inequality,
\begin{equation}\label{prop11proofb}
|E_{V_{\u}}H_1| + |E_{V_{\u}}H_2| = O\left\{\frac{1}{\sqrt{n}(1 + \eta^2)}\right\}.
\end{equation}
To bound $|E_{V_{\u}}H_3|$, we use integration by parts (Stein's lemma):
\begin{eqnarray*}
E_{V_{\u}}H_3 & = & \eta\frac{d}{n^{3/2}} E_{V_{\u}}\left\{\sum_{j = 1}^d \frac{\hat{\eta}^2s_j^{1/2}\tilde{\d}_j\tilde{u}_j}{(\hat{\eta}^2s_j + d/n)^2}\right\} \\
& = & \eta\frac{d}{n^{3/2}} E_{V_{\u}}\left\{\sum_{j = 1}^d \frac{\hat{\eta}^2s_j^{1/2}\tilde{\d}_j\tilde{u}_j}{(\hat{\eta}^2s_j + d/n)^2}\right\} \\
& = & 2\eta\frac{d}{n^{5/2}}E_{V_{\u}}\left\{\sum_{j = 1}^d \frac{(\hat{\eta}^2s_j - d/n)(\sqrt{n} \eta s_j^{1/2}\tilde{u}_j + \tilde{\d}_j)s_j^{1/2}\tilde{\d}_j\tilde{u}_j}{\hat{\s}^2(\hat{\eta}^2s_j + d/n)^3}\right\}.
\end{eqnarray*}
Thus,
\begin{eqnarray*}
|E_{V_{\u}}H_3| & \leq & 2\eta\frac{d}{n^{5/2}}E_{V_{\u}}\sum_{j = 1}^d \left|\frac{(\sqrt{n} \eta s_j^{1/2}\tilde{u}_j + \tilde{\d}_j)s_j^{1/2}\tilde{\d}_j\tilde{u}_j}{\hat{\s}^2(\hat{\eta}^2s_j + d/n)^2}\right| \\
& \leq & 8\eta\frac{d}{n^{5/2}}E_{V_{\u}}\sum_{j = 1}^d \left|\frac{(\sqrt{n} \eta s_j^{1/2}\tilde{u}_j + \tilde{\d}_j)\tilde{\d}_j\tilde{u}_j}{\hat{\s}^2(\hat{\eta}^2 + d/n)^2}\right|\left(\frac{1}{s_j^{3/2}} + s_j^{1/2}\right) \\
& \leq & 8 \eta^2\frac{d}{n^2}E_{V_{\u}}\left\{\frac{||\tilde{\dd}||}{\hat{\s}^2(\hat{\eta}^2 + d/n)^2}\left(\frac{1}{s_d} + s_1\right)\right\} \\
&& \qquad + 8\eta\frac{d}{n^{5/2}}E_{V_{\u}}\left\{\frac{\sqrt{\sum_{j = 1}^d \tilde{\d}_j^4}}{\hat{\s}^2(\hat{\eta}^2 + d/n)^2}\left(\frac{1}{s_d^{3/2}} + s_j^{1/2}\right)\right\} \\ 
& = & O\left\{\frac{1}{\sqrt{n}(1 + \eta^2)}\right\}.
\end{eqnarray*}
Combining this with (\ref{prop11proofa}) and (\ref{prop11proofb}) completes the proof of the proposition.  \hfill $\Box$

\vspace{.1in}

\textit{Proof of Proposition 11.}  Suppose that $0 < \th \leq d/n \leq \Th < 1$ for some fixed constants $\th,\Th \in \R$.  Let $\u \in \R^d$ be a unit vector and suppose that $\s^{-1}U\S^{1/2}\bb = \eta\u$, where $U$ is a $d \times d$ orthogonal matrix.  Define $\tilde{\S} = \tilde{\S}(\y,X) = U\S^{-1/2}\hat{\S}(\s\y,XU\S^{1/2})\S^{-1/2}U^T$ and let $\tilde{\bb}(\hat{\l}_r^*,\S) = \check{\bb}_r  = (X^TX + n \hat{\l}_r^*\S)^{-1}X^T\y$ be the adaptive ridge estimator defined in Proposition 10.  
Then 
\[
R_V\left\{\tilde{\bb}(\hat{\l}_r^*,\hat{\S})\right\} = R_{V_{\u}}\left\{\tilde{\bb}(\hat{\l}_r^*,\tilde{\S})\right\}
\]
and Proposition 10 implies that it suffices to show 
\[
R_{V_{\u}}\left\{\tilde{\bb}(\hat{\l}_r^*,\tilde{\S})\right\} = R_{V_{\u}}\left\{\tilde{\bb}(\hat{\l}_r^*,I)\right\} + O\left\{\left(E_{V_{\u}}||\tilde{\S} - I||^2\right)^{1/2}\right\}.
\]

Now notice that
\begin{eqnarray*}
\left|R_{V_{\u}}\left\{\tilde{\bb}(\hat{\l}_r^*,\tilde{\S})\right\} - R_{V_{\u}}\left\{\tilde{\bb}(\hat{\l}_r^*,I)\right\}\right|  & \leq & E_{V_{\u}}\left|||\tilde{\bb}(\hat{\l}_r^*,\tilde{\S})||^2 - ||\tilde{\bb}(\hat{\l}_r^*,I)||^2\right| \\
&& \quad + 2\eta E_{V_{\u}}\left|\left|\tilde{\bb}(\hat{\l}_r^*,\tilde{\S}) - \tilde{\bb}(\hat{\l}_r^*,I)\right|\right|.
\end{eqnarray*}
Considering integrands from the terms on the right-hand side above separately, we have
\begin{eqnarray}
\left|||\tilde{\bb}(\hat{\l}_r^*,\tilde{\S})||^2 - ||\tilde{\bb}(\hat{\l}_r^*,I)||^2\right|& = &\left|\frac{\hat{\eta}^4d^2}{n^4}\y^TX\left(\frac{\hat{\eta}^2}{n}X^TX + \frac{d}{n}I\right)^{-2}(I - \tilde{\S}^2) \right. \nonumber \\
&& \qquad \cdot \left(\frac{\hat{\eta}^2}{n}X^TX + \frac{d}{n}\tilde{\S}\right)^{-2}X^T\y \nonumber \\
&& \quad + \frac{\hat{\eta}^6d}{n^3}\y^TX\left(\frac{\hat{\eta}^2}{n}X^TX + \frac{d}{n}I\right)^{-2} \nonumber \\
&& \qquad \cdot \left\{\frac{1}{n}X^TX(I - \tilde{\S}) + (I - \tilde{\S})\frac{1}{n}X^TX\right\} \nonumber \\
&& \qquad \left.\cdot \left(\frac{\hat{\eta}^2}{n}X^TX + \frac{d}{n}\tilde{\S}\right)^{-2}X^T\y \right| \nonumber \\
& \leq & 4\left( s_1 +\frac{d}{n}\right) \frac{||X\hat{\bb}_{ols}||^2/n}{(\hat{\eta}^2 + 1)s_d^2}||I - \tilde{\S}|| \nonumber \\
& \leq &   \frac{4}{n-d}\left( s_1 + \frac{d}{n}\right) \frac{1}{s_d^2}||\y - X\hat{\bb}_{ols}||^2||I - \tilde{\S}|| \label{prop12proofa}
\end{eqnarray}
where $s_1 \geq \cdots \geq s_d$ are the eigenvalues of $n^{-1}X^TX$.  A similar calculation yields
\begin{equation}\label{prop12proofb}
||\tilde{\bb}(\hat{\l}_r^*,\tilde{\S}) - \tilde{\bb}(\hat{\l}_r^*,I)|| \leq \left(\frac{d}{s_dn} + 1\right)\frac{||X\hat{\bb}_{ols}||/n^{1/2}}{(\hat{\eta}^2 + 1)s_d^{1/2}}||I - \tilde{\S}||.
\end{equation}
The proposition follows by taking expectations in (\ref{prop12proofa})-(\ref{prop12proofb})  and applying Lemmas A2 and A3, along with the Cauchy-Schwarz inequality.  \hfill $\Box$

}\end{document}